\documentclass[preprint, 11pt, sort&compress]{elsarticle}
\usepackage[left=1in, right=1in, top=1in]{geometry}
\usepackage{hyperref}
\usepackage{mathrsfs}
\usepackage{nicefrac}
\usepackage{fancyhdr}
\usepackage{fancyhdr}
\usepackage{setspace}
\usepackage{lastpage}
\usepackage{upgreek}
\usepackage{graphicx}	
\usepackage[lofdepth,lotdepth]{subfig}
\usepackage{amsmath,amsthm}	
\usepackage{amssymb,dsfont,bbm}	
\usepackage{xcolor}  
\usepackage{algorithm2e}
\usepackage{bm}
\usepackage{xargs}
\usepackage[shortlabels]{enumitem}

\usepackage{aliascnt}
\usepackage{cleveref}

\usepackage{float}

\usepackage{comment}
\makeindex 

\setcitestyle{square}

\newcommandx{\AV}[2][1=]{{V}^{#1}_{\infty}(#2)}
\newcommandx{\SV}[2][1=]{{V}^{#1}_{n}(#2)}

\def\rmd{\mathrm{d}}
\def\rset{\mathbb{R}}
\def\ltwo{\mathrm{L}^2}

\numberwithin{equation}{section}

\theoremstyle{plain}

\newtheorem{Prop}{Proposition}
\newtheorem{Th}{Theorem}
\newtheorem{Lem}[Th]{Lemma}

\newtheorem*{prop:gen_conc}{Proposition~\ref{prop:gen_conc}}
\newtheorem*{prop:discr_conc}{Proposition~\ref{prop:discr_conc}}

\theoremstyle{remark} 

\newtheorem{As}{}
\newtheorem{Rem}{Remark}
\newtheorem{Corol}{Corollary}

\newenvironment{myassump}[2][]
{\renewcommand\theAs{\textbf{(#2)}}\begin{As}{#1}}{\end{As}}

\def\argmin{\operatornamewithlimits{argmin}}

\renewcommand\leq\leqslant
\renewcommand\geq\geqslant
\newcommand\eps\varepsilon

\newcommand{\bgamma}{{\boldsymbol{\gamma}}}

\newcommand{\R}{\mathbb R}

\newcommandx{\norm}[2][1=]{\ifthenelse{\equal{#1}{}}{\left\Vert #2 \right\Vert}{\left\Vert #2 \right\Vert^{#1}}}

\newcommandx{\stnorm}[2][1=]{\ifthenelse{\equal{#1}{}}{\biggl\Vert #2 \biggr\Vert}{\biggl\Vert #2 \biggr\Vert^{#1}}}

\newcommand{\AAA}{\mathsf{A}}
\newcommand{\BBB}{\mathsf{B}}
\newcommand{\fone}{f^{(1)}}
\newcommand{\ftwo}{f^{(2)}}
\newcommand{\Wone}{W^{(1)}}
\newcommand{\Wtwo}{W^{(2)}}
\newcommand{\vone}{v^{(1)}}
\newcommand{\vtwo}{v^{(2)}}
\newcommand{\requ}{\sigma}

\newcommand{\diagone}{\mathsf{\Delta}^{(1)}}
\newcommand{\diagtwo}{\mathsf{\Delta}^{(2)}}

\newcommand{\nset}{\mathbb{N}}
\newcommand{\nzeroset}{\mathbb{N}_0}

\newcommand{\Kxnorm}{\mathcal{K}}

\newcommand{\A}{\mathcal A}

\newcommand{\G}{\mathcal G}
\renewcommand{\H}{\mathcal H}

\newcommand{\NN}{\mathsf{NN}}

\def\sY{\mathsf{Y}}
\def\sZ{\mathsf{Z}}
\def\sT{\mathsf{T}}
\newcommand\logl[1]{\ell_{\text{log}}}
\newcommand\likel[1]{p_{\text{log}}}
\newcommand\pib[1]{\pi_{\text{log}}}
\newcommand\Ub[1]{U_{\text{log}}}
\newcommand{\fraca}[2]{#1/#2}

\def\sfp{\mathsf{p}}

\def\Cset{\mathsf{S}}

\def\Xset{\mathsf{X}}
\def\Xsigma{\mathcal{X}}

\def\ls{l}

\def\pg{P_{\ude}}
\def\hpow{\beta}
\def\mkp{\mathsf{P}}
\def\snzw{q}
\def\ude{\mathcal{D}}

\def\pdeg{\operatorname{deg}}
\def\udeg{\ude}

\def\epsnnapprox{\epsilon}

\def\hdeg{{\mathcal{D}_{H}}}

\newcommand{\indibr}[1]{\mathbbm{1}_{\{#1\}}}

\newcommand{\coint}[1]{\left[#1\right)}

\def\driftfunc{W}

\def\eqsp{\,}
\renewcommand{\H}{\mathcal{H}}

\def\sfp{\mathsf{p}}

\newcommandx{\PE}[2][1=]{\ensuremath{\mathsf{E}_{#1}\big[ #2 \big]}}

\newcommand{\E}{\mathsf{E}} 
\renewcommand{\P}{\mathbf{P}} 
\newcommand{\Vn}{V_{n}} 

\newcommandx{\hh}[1]{{\widehat{h}_{#1}}}

\newcommandx{\covcoeff}[3][1=\pi]{\rho_{#1}^{(#2)}(#3)}
\newcommandx{\ecovcoeff}[3][1=]{\hat{\rho}_{#1}^{(#2)}(#3)}

\newcommand{\ntest}{N}

\def\Var{\operatorname{Var}}

\def\PCov{\mathsf{cov}_\pi}
\def\tf{\tilde{f}}
\newcommand{\tvnorm}[1]{\| #1 \|_{\operatorname{TV}}}
\newcommandx\measureset[3][1=\mathrm{s},3=]{\mathbb{M}^{#3}_{#1}(\mathcal{#2})}

\newcommand{\Vnorm}[2]{\ensuremath{\left\Vert#1\right\Vert_{#2}}}
\newcommandx{\VnormFunc}[3][1=]{\ensuremath{\|#2\|_{{#3}}^{#1}}}
\def\driftfunc{W}

\newcommand{\lbeta}{\lfloor\beta\rfloor}

\def\rmd{\mathrm{d}}
\def\rme{\mathrm{e}}

\def\ltwo{L^2}

\usepackage{todonotes}

\def\eqsp{\,}
\newcommandx\PVar[1]{\operatorname{Var}_{{#1}}}

\def\Aset{{A_R}}

\def\metric{\mathsf{d}}

\newcommand{\1}{\ensuremath{\mathbbm{1}}}
\newcommandx{\indi}[2][1=]{\1^{#1}_{#2}}

\newcommandx\spacexd[3][1=\metric,3=]{\mathbb{S}_{#3}(\mathsf{#2},#1)}

\newcommand{\eqdef}{:=}
\newcommand{\setcomp}[1]{{({#1})^\mathsf{c}}}

\newcommand{\ggf}[2]{g_{{#1}, {#2} }}
\newcommand{\hhf}[2]{h_{{#1}, {#2} }}
\newcommand{\thhf}[2]{\tilde{h}_{{#1}, {#2} }}

\newcommand{\gcf}[2]{\mathcal{G}_{{#1}, {#2} }}
\newcommand{\hcf}[2]{\mathcal{H}_{{#1}, {#2} }}

\newcommand{\intg}[2]{\displaystyle\int\limits_{#1}^{#2}}
\newcommand{\smm}[2]{\displaystyle\sum\limits_{#1}^{#2}}

\DeclareMathAlphabet\mathbfcal{OMS}{cmsy}{b}{n} 

\crefname{theorem}{theorem}{Theorems}
\Crefname{Theorem}{Theorem}{Theorems}

\newaliascnt{lemma}{theorem}

\aliascntresetthe{lemma}
\crefname{lemma}{lemma}{lemmas}
\Crefname{Lemma}{Lemma}{Lemmas}

\newaliascnt{corollary}{theorem}

\aliascntresetthe{corollary}
\crefname{corollary}{corollary}{corollaries}
\Crefname{Corollary}{Corollary}{Corollaries}

\newaliascnt{proposition}{theorem}

\aliascntresetthe{proposition}
\crefname{proposition}{proposition}{propositions}
\Crefname{Proposition}{Proposition}{Propositions}

\newaliascnt{definition}{theorem}

\aliascntresetthe{definition}
\crefname{definition}{definition}{definitions}
\Crefname{Definition}{Definition}{Definitions}

\Crefname{As}{}{}
\crefname{As}{}{}

\newaliascnt{definitionProposition}{theorem}

\aliascntresetthe{definitionProposition}
\crefname{Proposition and Definition}{Proposition and Definition}{Proposition and Definition}
\Crefname{Proposition and Definition}{Proposition and Definition}{Proposition and Definition}

\newaliascnt{remark}{theorem}

\aliascntresetthe{remark}
\crefname{remark}{remark}{remarks}
\Crefname{Remark}{Remark}{Remarks}

\crefname{example}{example}{examples}
\Crefname{Example}{Example}{Examples}

\crefname{figure}{figure}{figures}
\Crefname{Figure}{Figure}{Figures}

\begin{document}


\begin{frontmatter}
\title{Theoretical guarantees for neural control variates in MCMC}
\author[inst1,inst2]{Denis Belomestny}
\ead{denis.belomestny@uni-due.de}
\address[inst1]{Duisburg-Essen University, Essen, Germany}
\author[inst2]{Artur Goldman}
\ead{agoldman@hse.ru}
\address[inst2]{HSE University, Moscow, Russia}
\address[inst3]{Institute for Information Transmission Problems, Moscow, Russia}
\author[inst2,inst3]{Alexey Naumov}
\ead{anaumov@hse.ru}
\author[inst2,inst3]{Sergey Samsonov}
\ead{svsamsonov@hse.ru}
\begin{abstract}
In this paper, we propose a variance reduction approach for Markov chains based on additive control variates and the minimization of an appropriate estimate for the asymptotic variance. We focus on the particular case when control variates are represented as deep neural networks. We derive the optimal convergence rate of the asymptotic variance under various ergodicity assumptions on the underlying Markov chain. The proposed approach relies upon recent results on the stochastic errors of variance reduction algorithms and function approximation theory.
\end{abstract}

\begin{keyword}
Variance reductions\sep MCMC\sep deep neural networks\sep control variates\sep Stein operator\,.
\end{keyword}
\end{frontmatter}

\section{Introduction}
\label{sec:introduction}
Modern statistical methods nowadays heavily use simulation-based numerical algorithms like Monte Carlo (MC) or Markov Chain Monte Carlo (MCMC). For example, Bayesian statistics relies upon  MCMC methods to calculate various functionals of the posterior distributions, see \cite{robert2004monte, rubinstein2016simulation}. Usually these simulation-based algorithms suffer from significant variance, and this is why variance reduction methods play an important role, see \cite{oates2017control, belomestny_empirical_2021, belomestny_variance_2020_esvm, oates2023}
\par
In many problems of computational statistics we aim to compute an expectation 
$\pi(f)\eqdef
\int_\Xset f(x)\pi(dx)
$,
where  $\pi$ is a distribution on $\Xset\subseteq\mathbb{R}^{d}$ and \(f:\Xset\to\mathbb{R}\), $f \in \ltwo(\pi)$.
Markov Chain Monte Carlo (MCMC) methods suggest to consider an estimate 
\begin{equation}
\label{eq:avg_def}
\textstyle
	\pi_n(f)\eqdef n^{-1}\sum_{k=0}^{n-1} f(X_k)
\end{equation}
based on a time homogeneous Markov chain \( (X_k)_{k \ge 0} \) taking values in \(\Xset\) with the corresponding Markov kernel $\mkp$ \cite{jones2004, roberts2004general, douc2018markov}. MCMC approach suggests to construct $\mkp$ in such a way that it admits the unique invariant distribution $\pi$. Under appropriate conditions, it can  be shown that for any $x \in \Xset$,
$\lim_{n \to \infty} \tvnorm{\mkp^n(x,\cdot) - \pi}=0$, where $\tvnorm{\cdot}$ is the total variation distance.
Moreover, the estimate \eqref{eq:avg_def} satisfies the central limit theorem (CLT)
    \begin{equation}
    \label{eq:clt_mcmc}
    \textstyle{
        \sqrt{n}\bigl[\pi_n(f) - \pi(f)\bigr]
        =
        n^{-1/2}
        \sum_{k=0}^{n-1}
        \bigl[f(X_k) - \pi(f)\bigr]\xrightarrow{d}
        \mathcal{N}(0, \AV{f})}
    \end{equation}
with the asymptotic variance of form
    \begin{equation}
    \label{eq:var_def_mcmc}
    \textstyle
    	\AV{f}=\lim_{n\to \infty} n\, \E_\pi\bigl[\{\pi_n(f) - \pi(f)\}^2\bigr]= 
    	\E_\pi\bigl[\tf^2\bigr] + 2 \sum_{k=1}^\infty \E_\pi\bigl[ \tf \mkp^k \tf\bigr]\,,
    \end{equation}
where we set $\tf = f-\pi(f)$. For conditions under which CLT holds, see
\cite{jones2004, roberts2004general, douc2018markov}.
Consequently, the true expectation $\pi(f)$ lies with probability $1-\delta$ in the asymptotic confidence interval
    \begin{equation}
        \label{eq:conf_interval}
    	\bigl( 
    	    \pi_n(f) - 
    	    c_{1-\delta/2}\sqrt{\frac{\AV{f}}{n}},
    	    \pi_n(f) +
    	    c_{1-\delta/2}\sqrt{\frac{\AV{f}}{n}}
    	\bigr),
    \end{equation}
where $c_{1-\delta/2}$ is $1-\delta/2$ quantile of standard normal distribution.
In order to make the above confidence interval tighter, we can either increase \(n\) or decrease \(\AV{f}.\) The latter problem can be solved by using variance reduction methods, in particular, the method of control variates \cite{mira2013zero, Assaraf1999, henderson1997variance, dellaportas2012control}. This method suggests constructing a class of control functions $\mathcal{G},$ that is, a class of functions \(g\) satisfying  $\pi(g)=0.$  Subtracting such functions do not introduce bias, that is $\pi(f)=\pi(f-g)$. Yet, variance of the new estimator $\pi_n(f-g)$ can be reduced by a proper  choice of $g$. Note that the closed form expression for \(V_\infty\) might be not available in many practical cases. Thus we can consider the optimisation problem 
\begin{equation}
\label{eq:optimisation_task}
\textstyle
\hat{g}_n \in \argmin_{g\in\G}\SV{f-g},
\end{equation}
where \(V_n\) is some estimate for the asymptotic variance \(V_\infty\). After solving \eqref{eq:optimisation_task} we consider a new estimate $\pi_N(f-\widehat{g}_n)$ for \(\pi(f)\) over new set of points $(X_0, \dots, X_{N-1})$. Note that $\AV{f-\widehat{g}_n}$ can be decomposed as

\begin{equation}
\label{eq:var_decomp}
\textstyle{
	\AV{f-\widehat{g}_n}
	=
	\underbrace{\AV{f-\widehat{g}_n}
	-
	\inf_{g\in\G}\AV{f-g}}_{\text{stochastic error}}
	+
	\underbrace{\inf_{g\in\G}\AV{f-g}}_{\text{approximation error}}}.
\end{equation}
Recently some works studied the  stochastic part of the error, see \cite{belomestny_variance_2020_esvm, belomestny_variance_2020_dependent, biz2018}. In this paper we focus both on stochastic and approximation errors in \eqref{eq:var_decomp} and obtain convergence rates of \(\AV{f-\widehat{g}_n}\) to zero as \(n\to \infty\) under a proper choice of the control variates class.
Let us stress that the term $\inf_{g\in\G}\AV{f-g}$ depends on the approximation properties of the class $\G$ with respect to the solution of the corresponding Poisson equation (see \Cref{sec: setup}). Although approximation error reduces when class \(\G\) gets larger, stochastic error grows while the complexity of the  class \(\G\) increases (see \Cref{th:main_slow}).
Thus, we can not simply reduce the overall error by taking some very large functional class \(\G\). In our paper, we show that an admissible choice of \(\G\) is given by the class of deep neural networks of a certain architecture.
\par
The main contribution of the paper is summarized in \Cref{th:bound_approx} and provides minimax optimal rates of variance reduction in terms of $V_{\infty}(f-g)$ for a Stein's control variate class based on deep neural networks. The key technical element of the proof is the relation between the term $\inf_{g\in\G}\AV{f-g}$ and the approximation properties of a given class of Stein control variates with respect to the solution of the corresponding Poisson equation.

The closest to our current study is the paper \cite{si_scalable_2020} where the authors performed a thorough study of popular control variates classes using the Stein identity. These classes are based on polynomials, kernel functions, and neural networks. The methodology suggested in \cite{si_scalable_2020} is based on the ideas developed in \cite{oates2017control, oates_convergence_2017} with additional empirical evaluations. The authors of \cite{si_scalable_2020} show   convergence  of the asymptotic variance for control variates based on polynomials and kernel functions. They also point out that proving convergence for  a neural network  based functional class is an open problem. 
Our technique is based on the approximation properties of certain classes of neural networks studied in
\cite{schmidt-hieber2020, belomestny2022simultaneous}. Approximation properties of neural networks will be the key to obtaining theoretic upper bound of asymptotic variance, which would not be possible while using functional classes of polynomials or kernel functions.
\par
The paper is organised as follows. In \Cref{sec: preliminaries} we describe used notation. In \Cref{sec: setup} we develop our approach. In \Cref{sec: assumptions} we formulate assumptions on which main theoretical results in \Cref{sec: theorres} are based.  Finally, in \Cref{sec:conclusion} we outline the significance of accomplished results and suggest improvements for future work. In \Cref{sec: proofres} we prove main theoretical results. Additional proofs are collected in \nameref{sec:appendix}. In particular, details on regularity of solutions to specific PDE's are provided in \ref{sec:pde_solution_regularity}. Section \ref{sec:excess_belomestny_statcomp} contains a statement of auxiliary results from \cite{belomestny_variance_2020_esvm}. Section \ref{sec:neural_networks_recu} contains the results on simultaneous approximation of a smooth function and its derivatives by neural networks with piecewise-polynomial functions. Additional details on numerical experiments are summarized in \ref{sec:numerical_exp_details}.

    \section{Notations and definitions}
    \label{sec: preliminaries}
    
    For simplicity we fix $\Xset = \rset^d$. If not stated otherwise, we denote sets as $\Xset = \rset^d$, $\Aset = (-R;R)^d$, $R>0$. Let $\setcomp{A}\eqdef\Xset\setminus A$, $\overline{A}$ is a closure of set $A$. For the sequences $(a_n)_{n \in \nset}$ and $(b_n)_{n \in \nset}$ we write $a_n \lesssim b_n$ if there exists a constant $C$ such that $a_n \leq C b_n$ for all $n \in \nset$. 
    \par 
    \paragraph*{\textbf{Norms and metrics}} Let \(\| \cdot \|\) denote the standard Euclidean norm. We say that \(f: \R^d \rightarrow \R\) is \(L-\)Lipschitz function if  \(|f(x)- f(x')| \le L \|x - x'\|\) for any \(x, x' \in \mathbb{R}^{d}\).
    \par
    Let $\driftfunc: \Xset \to \coint{1,\infty}$ be a measurable function.
    The $W$-norm of a function $h:\Xset\to\rset$ is defined as
    $\| h \|_{\driftfunc} = \sup_{x \in \Xset} \{|h(x)|/\driftfunc(x)\}$.
    For any two probability measures $\mu$ and $\nu$ on $(\Xset,\Xsigma)$ satisfying $\mu(\driftfunc) < \infty$ and $\nu(\driftfunc) < \infty$, the $W$-norm of $\mu-\nu$ is defined as $\Vnorm{\mu-\nu}{\driftfunc} = \sup_{\|f \|_{\driftfunc} \leq 1} |\mu(f) - \nu(f)|$.

    For a matrix $A$ and a vector $v$, we denote by $\|A\|_\infty$ and $\|v\|_\infty$  the maximal 
    absolute value of entries of $A$ and $v$, respectively.
    $\|A\|_0$ and $\|v\|_0$ shall stand for the number of non-zero entries of $A$ and $v$, 
    respectively. For a function $f : \Omega\to \mathbb{R}^d$ and a non-negative function $\sfp: \Omega\to \mathbb{R}$, we set 
    \begin{align*}
    	\|f\|_{L^\infty(\Omega)} = \sup_{x\in \Omega} \|f(x)\|,\ 
    	\|f\|_{L^2(\Omega)} =  \bigl(\int_{\Omega}\|f(x)\|^2\,\rmd x\bigr)^{1/2}, \|f\|_{L^p(\Omega, \sfp)} = \bigl( \int_{\Omega} \|f(x)\|^p \, \sfp(x) \, \rmd x \bigr)^{1/p}, p \geq 1.
    \end{align*}
    When the domain is clear, we omit $\Omega$ in the notations $L^\infty(\Omega)$, $L^p(\Omega)$, $L^p(\Omega, \sfp)$ and simply write $L^\infty$, $L^p$, and $L^p(\sfp)$, respectively, if there is no risk of ambiguity.
    
    \paragraph*{\textbf{Smoothness classes}}
    For any $s \in \nzeroset$, the function space $C^{s}(\Omega)$ consists of those 
    functions 
    over the domain $\Omega$ which have partial derivatives up to order $s$ in  $\Omega$, 
    and these derivatives are continuous in $\Omega$. If $0<\varpi \leq 1$, $u\in C^{s+\varpi}(\Omega)$ if $u\in C^s(\Omega)$ and derivative of order $s$ satisfy a H\"older condition on each compact subset of $\Omega$.
    For any multi-index $\bgamma = (\gamma_1,\dots,\gamma_d) \in 
    \nzeroset^d$, the partial differential operator $D^{\bgamma}$ is defined as
    \[
    	D^{\bgamma}f_i = \frac{\partial^{|\bgamma|} f_i}{\partial x_1^{\gamma_1} \cdots 
    	\partial x_d^{\gamma_d}}, \quad i \in \{1,\dots, m\}\,, \text{ and } \|D^{\bgamma} f\|_{L^\infty(\Omega)} = \max\limits_{1 \leq i \leq m} \|D^{\bgamma} f_i\|_{L^\infty(\Omega)}\,.
    \]
    \[
    \|f\|_{C^s} := \max\limits_{ 
    |\bgamma| \leq s} \|D^{\bgamma} f\|_{L^\infty(\Omega)}
    \]
    Here we have written $|\bgamma| = \sum_{i=1}^d \gamma_i$ for the order of 
    $D^{\bgamma}$.
    To avoid confusion between multi-indices and scalars, we reserve the bold font for the former ones.
    For the matrix of first derivatives, we use the usual notation $\nabla f = (\partial 
    f_i/\partial 
    x_j)$ \(i=1,\ldots,m\), $j = 1, \ldots, d$. For a function $\varphi: \R^d \mapsto \R$, $\varphi \in C^2(\Omega)$, we write $\Delta \varphi \in \R$ for its Laplacian. For a function $f:\Omega\to\R^m$ and any positive number $0 < \varpi \leq 1$, the
    \emph{H\"older constant} of order $\varpi$ is given by
    \begin{equation*}
    [f]_{\varpi}:=\max_{i \in \{1,\ldots,m\}}\sup_{x \not = y\in\Omega}\frac{|f_i(x)-f_i(y)|}{\min\{1, 
    \|x-y\|\}^{\varpi}}\;.
    \end{equation*} 
    Now, for any $\alpha >0$, we can define the \emph{H\"older ball} 
    $\H^{\alpha}(\Omega,H)$. If we let $s=\lfloor \alpha \rfloor$ 
    be the largest integer \emph{strictly less} than $\alpha$, $\H^{\alpha}(\Omega,H)$ contains functions in $C^s(\Omega)$ with $\varpi$-H\"older-continuous, $\varpi = \alpha - s > 0$, 
    partial derivatives of order $s$. Formally,
    \[
    \textstyle{
    	\H^{\alpha}(\Omega,H) = \big\{ f \in C^s(\Omega): \quad \|f\|_{\H^{\alpha}} := \max \{ 
    	\|f\|_{C^s}, \ \max\limits_{|\bgamma| = s} [D^{\bgamma}f]_{\varpi} \} \leq H \big\}.}
    \]
    Note, that due to $0<\varpi \leq 1$ norm $\|.\|_{\H^{s}}$ for $s\in\nset$ should be read as $\|f\|_{\H^{s}} = \|f\|_{\H^{(s-1)+1}}$, e.g. 
    
    $$\|f\|_{\H^{s}} = \max \{ 
    	\|f\|_{C^{s-1}}, \ \max\limits_{|\bgamma| = s-1} [D^{\bgamma}f]_{1} \}.
    	$$

Let $W^{k,p}(\rset^d)$ denote the Sobolev class of functions belonging to $L^p(\rset^d)$ along with their
generalized partial derivatives up to order $k$. Let $W^{k,p}_{loc}(\rset^d)$ denote the class of functions $f$ such that $\phi_0 \cdot f\in W^{k,p}(\rset^d)$ for all functions $\phi_0$ from the class $C_0^\infty(\rset^d)$ of infinitely differentiable functions with compact support.
\paragraph*{\textbf{Fixed point}} Given $\eps>0$, let  $\H_\eps\subset\H$ consist of centres of the minimal $\eps$-covering net of functional class $\H$ with respect
to the $\ltwo(\pi)$ distance. For a functional class $\mathcal{H}$ we define a fixed point as
\begin{equation}
\label{eq:definition-fixed-point}
\gamma_{\ltwo(\pi)}(\mathcal{H},n)
\eqdef
\inf \{ \eta>0: \
    		H_{\ltwo(\pi)}(\mathcal{H},\eta)\leq {n}\eta^2
\},
\end{equation}
where $H_{\ltwo(\pi)}(\mathcal{H},\eps) \eqdef \log |\H_\eps|$ and $|\H_\eps|$ is cardinality of $\H_\eps$.
    
\section{ESVM algorithm and control variates}
\label{sec: setup}
The key element of our approach to the asymptotic variance reduction is the optimisation problem \eqref{eq:optimisation_task}. To solve this problem we need to specify $V_n$, an estimate of $V_\infty$, and control variates class $\mathcal{G}$. Due to inherent serial correlation, estimating $V_\infty$ requires specific techniques such as spectral and batch mean methods; see \cite{flegal_galin} for a survey on variance estimators and their statistical properties. Following \cite{flegal_galin, belomestny_variance_2020_esvm}, we use the Empirical Spectral Variance estimator for $V_n$, which is defined below.
\par 
For a function $h: \Xset \mapsto \rset$ with $\pi(|h|)<\infty$ we set $\tilde{h} \eqdef h - \pi(h)$. For $s\in\nzeroset$ deﬁne the stationary lag $s$ autocovariance $\rho_\pi^{(h)}(s) \eqdef \E_\pi[\tilde{h}(X_0)\tilde{h}(X_s)] = \PCov(h(X_0), h(X_s))$ and the lag $s$ sample autocovariance via
\begin{equation}
\label{eq:empirical-autorcovariance}
\textstyle{
\ecovcoeff[n]{h}{s} \eqdef n^{-1} \sum_{k=0}^{n-s-1} \{h(X_k) - \pi_n(h)\} \{h(X_{k+s}) - \pi_n(h)\}},
\end{equation}
where $\pi_n(h) \eqdef  n^{-1} \sum_{j=0}^{n-1} h(X_j)$. The spectral variance estimator is based on  truncation and weighting of the sample autocovariance function,
\begin{equation}
\label{eq:sv}
\textstyle{
\SV{h} \eqdef \sum_{s=-(b_n-1)}^{b_n-1} w_n(s) \ecovcoeff[n]{h}{|s|}}\eqsp,
\end{equation}
where $w_n(s)$ is the \emph{lag window} and $b_n\in\nset$ is the \emph{truncation point}. Here
the lag window is a kernel of the form \(w_n(s)=w(s/b_n)\), where \(w(\cdot)\) is a symmetric non-negative function supported on \([-1,1]\) which satisfies $|w(s)| \leq 1$ for \(s \in [-1,1]\) and \(w(s)=1\) for \(s \in [-1/2,1/2]\).
\par 
The described estimator is used in the Empirical Spectral Variance Minimisation algorithm (ESVM), introduced in \cite{belomestny_variance_2020_esvm}. The ESVM algorithm with general class of control variates $\G$ is summarised in \Cref{algorithm:esvm}:
\par 
\RestyleAlgo{ruled}
\begin{algorithm}[H]
\setstretch{1.35}
\caption{ESVM method}
\label{algorithm:esvm}
\KwIn{Two independent sequences: $\mathbf{X}_n=(X_k)_{k=0}^{n-1}$ - for train, $\mathbf{X}'_{\ntest}=(X'_k)_{k=0}^{\ntest-1}$ for test;}
\textbf{1.} Choose a class $\G$ of functions with $\pi(g)=0$ for all functions $g\in\G$\;
\textbf{2.} Find  $\widehat{g}_n\in\argmin_{g\in\G}\Vn(f-g)$, where $\Vn(\cdot)$ is computed according to \eqref{eq:sv}\;
\KwOut{\(\pi_{\ntest}\left(f - \widehat{g}_n\right) \) 
computed using $\mathbf{X}'_{\ntest}$.}
\end{algorithm}
Paper \cite{belomestny_variance_2020_esvm} analyses stochastic error for the general class $\G$. The aim of our analysis is to study the properties of the approximation error in \eqref{eq:var_decomp} based on the expressive power of the class $\G$. Moreover, for the particular case of neural networks-based control variates we show, how to choose the network's architecture in order to balance the \emph{stochastic} and \emph{approximation} errors in \eqref{eq:var_decomp}. Our analysis yields minimax optimal variance reduction rates, see \Cref{th:bound_approx}.
\par 
In \Cref{algorithm:esvm} we need to fix a class of control variates $\mathcal{G}$. Following  \cite{Assaraf1999, mira2013zero, oates_convergence_2017} we choose $\mathcal{G}$ as a class of Stein control variates of the form 
\begin{equation}
\label{eq:steincv}
\textstyle{
g_\varphi(x) \eqdef \Delta \varphi(x) +
\langle \nabla\log \pi(x), \nabla \varphi(x) \rangle}\eqsp,
\end{equation}
    where $\varphi \in C^2(\Xset)$.
    In practice it would be convenient to restrict the choice of functions $\varphi$ to some subclass of $C^2(\Xset)$. In particular, some recent works, which have studied a practical performance of the control variates method \cite{belomestny_variance_2020_esvm, si_scalable_2020, south2023regularized} have used classes of polynomials, kernel functions and neural networks.
    \par
    Under rather mild conditions on $\pi$ and $\varphi$, it follows from integration by parts that $\pi(g_\varphi)=0$; see \cite[Propositions~1 and 2]{mira2013zero}. 
    If there exists $\varphi^*$ which satisfies the relation
    \begin{equation}
    \label{eq:phi_star_def}
        g_{\varphi^*}(x) =
        f(x)-\pi(f),
    \end{equation}
    then $g_{\varphi^*}$ is the optimal control variate in a sense that $\AV{f-g_{\varphi^*}} = 0$.

    However, in practice we can not compute $g_{\varphi^*}$. If $\pi(f)$ was available, and we were able to evaluate $f(x)$ in every point $x$, we would have an access to $g_{\varphi^*}(x)$ in every point. Unfortunately this is not the case, because out initial aim is to estimate $\pi(f)$ itself, and thus having access to it is unrealistic. Moreover, in practice we might not have an access to an analytic expression form of $f(x)$, but rather we could evaluate $f(x)$ only for a fixed set of points $x$. Access to $g_{\varphi^*}$ through $\varphi^*$ is even more unrealistic, because for this we would not only need to know $\pi(f)$ and $f(x)$, but  also to solve the elliptic differential equation \eqref{eq:phi_star_def}. In fact, the solution of this equation is not available in the closed form and can only be approximated. Both  analysis of the properties of $\varphi^*$ and its approximation by  parametric families of functions are non-trivial problems taking into account  non-boundedness of the set $\Xset$. For this reason, we suggest to introduce a bounded set $\Aset = (-R;R)^d$ and consider the equation
    \begin{equation}
    \label{eq:diff_eq_alt}
        \Delta \varphi^*_{R}(x)
        +
        \langle \nabla\log \pi(x), \nabla \varphi^*_{R}(x) \rangle
        =
        f(x)-\pi(f)\,, x \in \Aset.
    \end{equation}
    For further theoretical analysis we introduce and use a class of \emph{biased} control variates
    \begin{equation}
    \label{eq:ggf}
        \ggf{F}{R}(x) = 
        \begin{cases}
        \Delta F(x)
        +
        \langle \nabla\log \pi(x), \nabla F(x) \rangle
        ,\ x\in \Aset,\\
        0,\ x\in \setcomp{\Aset}\,.
        \end{cases}
    \end{equation}
    for any $F \in C^2(\Xset)$. Note that $\ggf{F}{R}(x) = g_F(x)\indibr{x \in \Aset}$, $g_F$ is as defined in \eqref{eq:steincv}. Similarly,
    \begin{equation}
    \label{eq:hhf}
        \hhf{F}{R}(x) \eqdef f(x)-\ggf{F}{R}(x)\,,
    \end{equation}
    and for an arbitrary function class $\Phi,$ we introduce
    \begin{equation}
        \label{eq:f_minus_steincv}
        \hcf{\Phi}{R} = \{\hhf{\varphi}{R} = f-\ggf{\varphi}{R}\ :\ \ggf{\varphi}{R}\in\gcf{\Phi}{R}\},\ 
        \gcf{\Phi}{R} = \{\ggf{\varphi}{R}: \varphi\in \Phi\}\,.
    \end{equation}
    When there is no risk of confusion we will omit the replicating index $R$ and write $g_{\varphi^*_{R}}$ instead of $g_{\varphi^*_{R}, R}$. Note that $g_{\varphi^*_{R}}$ is no longer optimal control variate, that is, $\AV{f-g_{\varphi^*_{R}}} \neq 0$. Moreover, in general, $\pi(g_{\varphi^*_{R}}) \neq 0$. Yet we can account for the  bias $\pi(g_{\varphi^*_{R}})$ by a proper choice  of the parameter $R$.
       
    \section{Assumptions}
     \label{sec: assumptions}
    In this section, we formulate assumptions on the considered Markov Chain $(X_k)_{k \geq 0}$, the target distribution $\pi$ and the target function $f$. We consider the target distribution $\pi$ which has a bounded density w.r.t. the Lebesgue measure given by
    \begin{equation}
    \label{eq:pi_def}
    \textstyle{
    \pi(x)= e^{-U(x)}\Big/\int_\Xset e^{-U(y)} dy}\eqsp.
    \end{equation}
    In terms of the introduced potential $U(x),$ the equation \eqref{eq:diff_eq_alt} writes as 
    \begin{equation}
        \Delta \varphi^*_{R}(x)
        -
        \langle \nabla U(x), \nabla \varphi^*_{R}(x) \rangle
        =
        f(x)-\pi(f), \quad x \in \Aset.
    \end{equation}
    Thus, the smoothness properties of $\varphi^*_{R}$ depend upon the smoothness of $f$ and $U$. 
    In particular, we need to control the growth rate of $\nabla U$ and $f$. To this end, we impose the following assumptions.     
    \begin{myassump}{AUF}
    \label{assu:AUF}
    The potential $U$ is $\mu$-strongly convex and the function $f$ is bounded with $\norm{f}_{L^\infty(\Xset)}=B_f<\infty$.
    Moreover, there exist $C_0, R_0, \ude \geq 0$, $0<\hpow<\infty$ such that
    \[
    f, \nabla U \in \mathcal{H}^{\hpow+1}(\overline{\Aset}, \pg(R))
    \]
    for all $R>0$ where for $R\geq R_0,$ it holds that $\pg(R)\leq C_0 \cdot R^\ude$.
    \end{myassump}
    In the assumptions above, $\mu$-strong convexity of $U$ is required to control the bias of control variates $g_{F,R}$ introduced in \eqref{eq:ggf}. Smoothness constraints on $U$ and $f$ are introduced in order to control the growth rate and the smoothness of $\varphi^*_{R}$. The requirement that  $f$ is bounded does not appear in our bound for the approximation error in \eqref{eq:var_decomp}, but it is essential to control the corresponding stochastic error using \Cref{th:main_slow}.
    \paragraph{\textbf{Markov Chain assumptions}}
    Now we state assumptions on the Markov kernel $\mkp$, which are \emph{geometric ergodicity} \Cref{assu:ge} and \emph{bounded recurrence} \Cref{assu:br}. Let $W:\Xset\to[1;\infty)$ be a measurable function.
    
    \begin{myassump}{GE}
    \label{assu:ge}
    The Markov kernel $\mkp$  admits a unique invariant probability measure \(\pi\) such that $\pi(W)<\infty$ and there exist  $\varsigma >0, 0< \rho < 1$ such that  for all $x \in \Xset$ and $n \in \nset,$
    \begin{equation*}
    \Vnorm{\mkp^n(x,\cdot)-\pi}{\driftfunc}\leq \varsigma \driftfunc(x) \rho^{ n}.
    \end{equation*}
    \end{myassump}

 \begin{myassump}{BR}
    \label{assu:br}
    There exist a non-empty set $\Cset \subset \Xset$ and real numbers $u>1, J>0$ and $\ls>0$ such that
    \begin{equation*}
    \sup_{x \in \Cset} \E_x[u^{-\lambda}] \le J \quad \text{ and } \quad \sup_{x \in \Cset} W(x) \le \ls,
    \end{equation*}
    where $\lambda$ is the return time to the set $\Cset$.
    \end{myassump}
    Considering the $W$-geometrically ergodic Markov kernels is rather standard in the MCMC literature, see e.g. \cite{adamczak2015exponential, lemanczyk2020general}. The main advantage of this type of ergodicity is that it allows to write quantitative deviation bounds for additive functionals of Markov chains, see \cite[Section~23]{douc2018markov}.

    \section{VR rates with ESVM procedure}
    \label{sec: theorres}
    Before we prove the variance bound, we must specify the underlying class $\Phi$ in \eqref{eq:f_minus_steincv}. In this study, we consider $\Phi$ as a particular class of fully-connected deep neural networks. Note that there is a number of recent studies devoted to the approximation of functions by neural networks with sufficiently smooth activation function, see \cite{guhring_approximation_2020, de2021approximation, langer2021approximating}. There are also recent works studying the problem of approximating PDE solutions in  a suitable norm, see \cite{han2018solving, grohs2022deep}. 
    \par 
    Our considered application has some specific requirements. First, we need to control the statistical complexity of the underlying class of neural networks. In particular, this requires the weights of the neural network to be bounded (say by $1$), similar to \cite{schmidt-hieber2020}. Second important question is the choice of activation function. The approximation network should be substituted into the Stein operator \eqref{eq:steincv}, which imposes smoothness constraints. Thus, we need to have at least two times continuously differentiable functions $\varphi$, otherwise, we are not guaranteed to have $\pi(g_\varphi)=0$. As a result, we choose a piecewise-polynomial activation function, namely, the rectified cubic unit $\sigma^{\mathsf{ReCU}}(x) = (x \vee 0)^3$ . The approximation properties of networks with piecewise-polynomial activations are well-studied, see  \cite{powernet2020, li2022, siegelReluk, belomestny2022simultaneous}. The main result of this section will rely on the appropriate modification of \cite{belomestny2022simultaneous} and careful control of the covering number of the respective class of networks. 
    \par 
    Now we formally define the underlying class of neural networks. We first fix an activation function 
    $\sigma: \R \rightarrow \R$. For a  vector $v= (v_1,\dots,v_p) \in \R^p$,
    we define the shifted activation function $\sigma_v: \R^p \rightarrow \R^p$ as 
    \begin{equation*}
    	\sigma_{v}(x) = \bigl(\sigma(x_1-v_1),\dots,\sigma(x_p-v_p)\bigr), \quad x = 
    	(x_1,\dots,x_p) 
    	\in 
    	\R^p.
    \end{equation*}
    Given a positive integer $L$ and a vector $\A = (p_0, p_1, \dots, p_{L+1}) \in \mathbb N^{L+2}$,  a neural network of depth $L+1$ (with $L$ hidden layers) and architecture \(\A\) is  a function of the form
    \begin{equation}
    	\label{eq:nn}
    	f: \R^{p_0} \rightarrow \R^{p_{L+1}}\,, \quad 
    	f(x) = W_L \circ \sigma_{v_{L}} \circ W_{L-1} \circ \sigma_{v_{L-1}} \circ \dots \circ W_1\circ\sigma_{v_1} \circ 
    	W_0 \circ x,
    \end{equation}
    where $W_i \in \R^{p_{i+1} \times p_i}$ are weight matrices and $v_i \in \R^{p_i}$ are shift vectors. Next, we introduce a special subclass class of neural networks of depth $L+1$ with architecture $\A$:
    \[
    	\NN(L, \A) =
    	\left\{ \text{$f$ of the form \eqref{eq:nn}} : \|W_0\|_\infty \vee \max\limits_{1 \leq \ell 
    	\leq L} \left\{ \|W_\ell\|_\infty \vee \|v_\ell\|_\infty \right\} \leq 1\right\}.
    \]
    The maximum number of neurons in one layer $\|\A\|_\infty$ is called the width of the neural network.
    We use the ReCU (rectified cubic unit) activation function, defined as
    \begin{equation}
            \label{eq:recu_def}
    	\sigma^{\mathsf{ReCU}}(x) = (x \vee 0)^3.
    \end{equation}
    From now on, we use ReCU activation functions and write $\sigma(x)$ instead of $\sigma^{\mathsf{ReCU}}(x)$ for short. However, further in the proof the following functions will appear:
    \begin{equation}
    \label{eq:requ_relu_def}
    \sigma^{\mathsf{ReQU}}(x) = (x \vee 0)^2,\quad
    \sigma^{\mathsf{ReLU}}(x) = (x \vee 0)
    \end{equation}  
    In what follows, we consider sparse neural networks assuming that only a few weights are not equal to zero. For this purpose, we introduce a class of neural networks of depth $L+1$ with architecture $\A$ and with at most $\snzw$ non-zero weights:
    \begin{equation}
        \label{eq:nn_cl_def}
        \textstyle{
    	\NN(L, \A, \snzw) =
    	\bigl\{ f \in \NN(L, \A) : \|W_0\|_0 + \sum_{\ell = 1}^L \left(\|W_\ell\|_0 + 
    	\|v_\ell\|_0 \right) \leq \snzw\bigr\}}\eqsp.
    \end{equation}
    \par 
    Unfortunately, the class \eqref{eq:nn_cl_def} might be overly too complex to be used as a basic class $\Phi$ for Stein control variates $\gcf{\Phi}{R}$ defined in \eqref{eq:f_minus_steincv}. It is because the architecture $\A$ may have too many non-zero weights in order to ensure good approximation properties in $\|\cdot\|_{\H^{3}(\Aset)}$ for large values of $R$. In particular, without additional constraints, the class \eqref{eq:nn_cl_def} will contain functions $f(x)$ which grow too fast for large values of $\|x\|$. 
    \par 
    For this reason we consider as an underlying class $\Phi_R$ the intersection of $\NN(L, \A, \snzw)$ with a ball of suitable radius in $\|\cdot\|_{\H^{3}(\Aset)}$ norm. More precisely, in order to define $\Phi_R$, we introduce for parameters $R,K > 0$ and absolute constant $C_{\ref{prop:recu_approx_const}, 1}$ defined in \Cref{prop:recu_approx_const}, the quantities:
    \begin{equation}
    \label{eq:eps_nnaprox_def}
    H_{R} =  \|\varphi^*_R\|_{\H^{(\hpow+3)}(\Aset)},\
    \epsnnapprox = \frac{C_{\ref{prop:recu_approx_const}, 1}^{d(\hpow+3)}(2R)^{(\hpow+3)} H_{R} (\hpow+3)^3 }{K^{\hpow}}.
    \end{equation}
    Then, with $v
    =
    \frac{1}{2}(1,1,\dots,1)^T$, we define 
    \begin{equation}
    \label{eq:phi_def}
    \Phi_{R}
    =
    \biggl\{
    \varphi = g\circ
    \left(
    \frac{x}{2R}+v
    \right)
    \biggm| \norm{
    \varphi
    }_{\H^3(\Aset)}\leq 2((2R)^{\hpow+3} H_{R} + \epsnnapprox
    ), \nonumber
     g\in \NN (L, \A(H_{R}, K),\snzw(H_{R}, K))
    \biggr\},
    \end{equation} 
    where the precise values of $L,\ \A(H_{R}, K),\ \snzw(H_{R}, K)$ are outlined in \Cref{prop:recu_approx_const}. Properties of the constructed class are discussed in \Cref{th:nn_constr}. From now on, we use ESVM algorithm with the introduced class $\Phi_R$. For the respective class of control variates $\hcf{\Phi_R}{R}$ defined in \eqref{eq:f_minus_steincv}, we omit the repeating index $R$ and write simply $\mathcal{H}_{\Phi_R}$.
    \par 
    For our theoretical analysis, instead of looking for a function with the smallest spectral variance in the whole class $\mathcal{H}_{\Phi_R}$ we will perform optimization over a finite $\eps$-covering net in $\mathcal{H}_{\Phi_R}$, denoted as $\mathcal{H}_{\Phi_R,\eps}$.
    Epsilon-nets are introduced to simplify our theoretical analysis and  both control variates (picked  in the whole class $\H$ and in its $\eps$-net) have similar properties in general. Of course, if class $\H$ is parametric, there is no need to ``discretize" it in order to perform  \Cref{algorithm:esvm}. It is more natural to use some gradient-based optimization algorithm to find  minimum over $\H$.
    Nevertheless, if $\H$ is non-parametric, 
    numerical optimization over the whole class $\H$ is generally impossible. For instance, in \cite{biz2018} we have a lengthy discussion motivating the use of $\eps$-nets in Section 2. There we also present a concrete example of constructing  an $\eps$-net based on neural nets.  This type of nets can be directly used in practice. 
    In the sequel, we assume that $\mathcal{H}_{\Phi_{R},\eps}$ is finite. Then, instead of solving the optimisation problem \eqref{eq:optimisation_task}, we consider its relaxed version
    \begin{align}
    \label{eq:vmnq}
    \widehat{h}_{n,\eps, R} \in \argmin_{h \in \mathcal{H}_{\Phi_R,\eps}} \SV{h}\,.
    \end{align}
    With the notation above, we introduce the corresponding NN-based theoretical procedure, which is further referred to as ESVM-NN. 
    \par 
    \RestyleAlgo{ruled}
    \begin{algorithm}[H]
    \setstretch{1.35}
    \caption{Algorithm ESVM-NN} 
    \label{algorithm:esvm_nn}
    \KwIn{Two independent Markov chains with same kernel $\mkp$, satisfying \Cref{assu:ge} and \Cref{assu:br}: $\mathbf{X}_n=(X_k)_{k=0}^{n-1}$ - for train, $\mathbf{X}'_{\ntest}=(X'_k)_{k=0}^{\ntest-1}$ for test;}
    \textbf{1.} Set $R=\log n$, $K = n^{\frac{1}{2\hpow + d}}$, \(b_n=2(\log(1/\rho))^{-1}\log(n)\), \\ 
    \textbf{2.} Choose a class $\G=\gcf{\Phi_R}{R}$ and construct the corresponding class $\mathcal{H}_{\Phi_R}$\;
    \textbf{3.} Set $\eps = \gamma_{\ltwo(\pi)}(\mathcal{H}_{\Phi_R},n)$\;
    \textbf{4.} Find  $\widehat{h}_{n,\eps, R} \in \argmin_{h \in \mathcal{H}_{\Phi_R,\eps}} \SV{h}$ where $\Vn(\cdot)$ is computed according to \eqref{eq:sv}\;
    \KwOut{\(\pi_{\ntest}(\widehat{h}_{n,\eps, R}) \) 
        computed using $\mathbf{X}'_{\ntest}$.}
    \end{algorithm}
    
    Now we solve \eqref{eq:vmnq} for the biased control variates from \eqref{eq:ggf}, \eqref{eq:hhf} using the class $\Phi_{R}$ defined in \eqref{eq:phi_def}. The next theorem is the main result of our study providing convergence rates of the estimate \(\widehat{h}_{n,\eps, R}\) in terms of its asymptotic variance.
    \begin{Th}
    \label{th:bound_approx} 
    Assume \Cref{assu:AUF}, \Cref{assu:ge}, and \Cref{assu:br}. Then for any $x_0 \in \Cset$ and $\delta \in (0,1)$ there exists \(n_0=n_0(\delta, R_0, d, \hpow, \udeg)>0\) such that for all $n\geq n_0$, $n\in\nset$ by setting $R=\log n$, $K = n^{\frac{1}{2\hpow + d}}$, \(b_n=2(\log(1/\rho))^{-1}\log(n)\), $\eps = \gamma_{\ltwo(\pi)}(\mathcal{H}_{\Phi_R},n)$ with $\P_{x_0}-$probability at least $1-\delta$, it holds that 
    \[
        \AV{\widehat{h}_{n,\eps, R}}
        \lesssim
        (\log n)^{C_{\ref{th:bound_approx}, 1}} \max\left(n^{-\frac{\hpow}{2\hpow+d}}, \frac{\log(\fraca{1}{\delta})}{\sqrt{n}}\right)\eqsp,
    \]
    where $C_{\ref{th:bound_approx}, 1}=C_{\ref{th:bound_approx}, 1}(\hpow, \udeg)$ is independent of the problem dimension $d$. Moreover, with probability at least $1-\delta$ it holds
    \[
    \left|
    \pi(f)-\pi_N(\widehat{h}_{n,\eps, R})
    \right|\lesssim
    c_{1-\fraca{\delta}{2} } \sqrt{\frac{\AV{\widehat{h}_{n,\eps, R}}}{N}}+ (\log n)^{C_{\ref{th:bound_approx}, 2}}\cdot n^{C_{\ref{th:bound_approx}, 3}-\log n},
    \]
    where $c_{1-\fraca{\delta}{2}}$ is $1-\delta/2$ quantile of standard normal distribution and $C_{\ref{th:bound_approx}, 2} = C_{\ref{th:bound_approx}, 2}(\hpow, d, \udeg)$, $C_{\ref{th:bound_approx}, 3}=C_{\ref{th:bound_approx}, 3}(\hpow, d)$. Here $\lesssim$ stands for inequality which holds up to multiplicative constants, depending on $(\varsigma, W, x_0, J, l, \rho, u, C_0, \hpow, d, B_f, \mu)$.
    \end{Th}

\section{Numerical study}
\label{sec: numstudy}
In this section, we study the numerical performance of the described method for various distributions. PyTorch based Python implementation is available at \href{https://github.com/ArturGoldman/NNESVM}{https://github.com/ArturGoldman/NNESVM}. For each experiment we need to fix the following parameters: 
\begin{itemize}
\item kernel function $w(s)$
\item lag-window size $b_n$
\item size of training trajectory $n=n_{\text{burn}}+n_{\text{train}}$
\item size of testing trajectory $n_{\text{test}}$ and number of testing trajectories $T$
\item chain generation method (CGM), its parameters 
\item dimensionality of data $d$
\item function $f$ of interest
\end{itemize}
	
	For CGM we used Unadjusted Langevin Algorithm (ULA) and No-U-Turn Sampler (NUTS) which are described in \cite{belomestny_variance_2020_esvm} and \cite{hoffman_no-u-turn_2011}, respectively.
	
	If not specified, by default we fix triangular function as a kernel function $w(s)$ supported on $[-1;1]$

	\begin{equation*}
	    w(s)
	    =
	    \begin{cases}
	        1+s,\ -1\leq s < 0,\\
	        1-s,\ 0 \leq s \leq 1
	    \end{cases}
	\end{equation*}
	
	Also we fix $f(X) = X_2^2$ as a function of interest and calculate $\E_\pi{X_2^2}$ for synthetic datasets. For Logistic Regression, we evaluate average likelihood over the test set. For target distribution $\pi,$ we define potential function $U$ as in \eqref{eq:pi_def}. All other necessary hyperparameters for all experiments are reported in \ref{sec:numerical_exp_details}.
    We report Empirical Spectral Variance reduction rate $V_n(f)/V_n(f-g)$ (ESVRR) and boxplots  representing empirical variance on test chains.
    In the following we write ReLU, ReQU and ReCU for the activation functions $\sigma^{\mathsf{ReLU}}(x)$, $\sigma^{\mathsf{ReQU}}(x)$, $\sigma^{\mathsf{ReCU}}(x)$ respectively, which are defined in \eqref{eq:recu_def}, \eqref{eq:requ_relu_def}.
	
	\subsection{Implementation specifics}
	Before delving into the experiments, let's outline the technical details and implementation specifics. Initially, we employ a one-layered fully-connected neural network (MLP) with the Rectified Cubic Unit (ReCU) activation function. In certain implementations of Stein control variates, an alternative neural network structure outputs \(\nabla \varphi\) instead of \(\varphi\), as seen in works such as \cite{si_scalable_2020, oates2023}. This design choice eliminates the need for calculating second-order derivatives, requiring only a single backpropagation to evaluate \(\Delta \varphi\). This approach accelerates training and conserves memory, given the necessity of differentiably evaluating \(\Delta \varphi\). However, in this context, we opt for a neural network that approximates the function \(\varphi\), necessitating an additional backpropagation procedure for updating neural network weights. Consequently, even relatively straightforward experiments entail substantial computational demands, emphasizing the computational complexity even with basic neural network architectures. This experimental setup aligns with our goal of ensuring consistency with the theoretical framework developed earlier.
	\par
	In our experiments, we employ the Adam optimizer with a small exponential weight decay. It's essential to note that our study does not delve into identifying the optimal optimization strategy that yields the maximum possible variance reduction. Additionally, we do not extensively explore regularization techniques that could be incorporated into the loss function to enhance the stability of the learning process (for more details, refer to \cite{si_scalable_2020}). Further information on the optimization methods utilized can be found in the implementation details, emphasizing the focus of this work on the core aspects of variance reduction and neural network approximation rather than optimization strategy fine-tuning.
	\par
	A notable implementation detail involves the generation of the Markov chain, which typically necessitates a consecutive construction. However, we introduce the possibility of parallelizing the chain generation process and the subsequent evaluation of statistics. Additionally, we've incorporated the option to utilize pregenerated Markov chains. This feature aims to save time during evaluation epochs by eliminating the need for on-the-fly chain generation, providing flexibility in the experimental setup.
	
	\subsection{Funnel distribution}
	
	Here we study Funnel distribution from \cite{neal_slice_2000}. It has density the form
	
	\begin{equation*}
	    p(x_1, \dots, x_n)
	    =
	    p_\mathcal{N}(x_1|0, a)
	    \prod_{i=2}^d
	    p_\mathcal{N}(x_i|0, e^{2bx_1}), \quad a,b>0
	\end{equation*}
	(here $p_\mathcal{N}(x|\mu,\sigma^2)$ denotes the density of Normal distribution with mean $\mu$ and variance $\sigma^2$) and its potential function can be written as
	
	\begin{equation*}
	    U(x_1, \dots, x_n)
	    =
	    \frac{x_1^2}{2a}
	    +
	    (d-1)bx_1
	    +
	    e^{-2bx_1}\cdot 
	    \smm{i=2}{d}
	    \frac{x_i^2}{2}
	\end{equation*}
    We take $d=2$, $a=1$,
        $b=0.5$ and we compare the performance of different methods. In particular, we compare neural networks with different activation functions. Also, we test linear regression on polynomials  with degree from $0$ up to $4$. All models were trained with similar optimisation procedures and the the same hyperparameters and  number of steps. The hyperparameters are shown in Table~\ref{tab:funnel}

        \begin{table}[H]
        \centering
        \begin{tabular}{ |c|c|c|c|c|c|c| } 
        \hline
        CGM & 
        $\gamma$&
        $n_{\text{burn}}$& 
        $n_{\text{train}}$&
        $n_{\text{test}}$&
        $T$&
        $b_n$\\
        \hline
        NUTS & 
        0.1&
        $10^4$& 
        $3\cdot 10^4$&
        $3\cdot 10^4$&
        30&
        30\\
        \hline
        \end{tabular}
        \caption{\label{tab:funnel} Funnel, algorithm and chain generation hyperparameters}
    \end{table}
    \par
    It is worth mentioning that if in the same setting neural network with ReLU activation function is used, a noticeable bias appears. We believe that this bias is due to the network not being  2 times differentiable.
	
	\begin{table}[H]
 \centering
    \begin{tabular}{|c|c|}
    \hline
        Setup & $V_n(f)/V_n(f-g)$  \\ 
    \hline
    Linear regression, degree 4 & 4.9\\
    \hline
     ReQU activation & 4.8\\
    \hline
      ReCU activation & 15.9\\
    \hline
    \end{tabular}
    \caption{Funnel, results}
    \end{table}

    \begin{figure}[H]
      \centering
      \subfloat[Boxplots]
      {
        \includegraphics[width=0.48\textwidth]{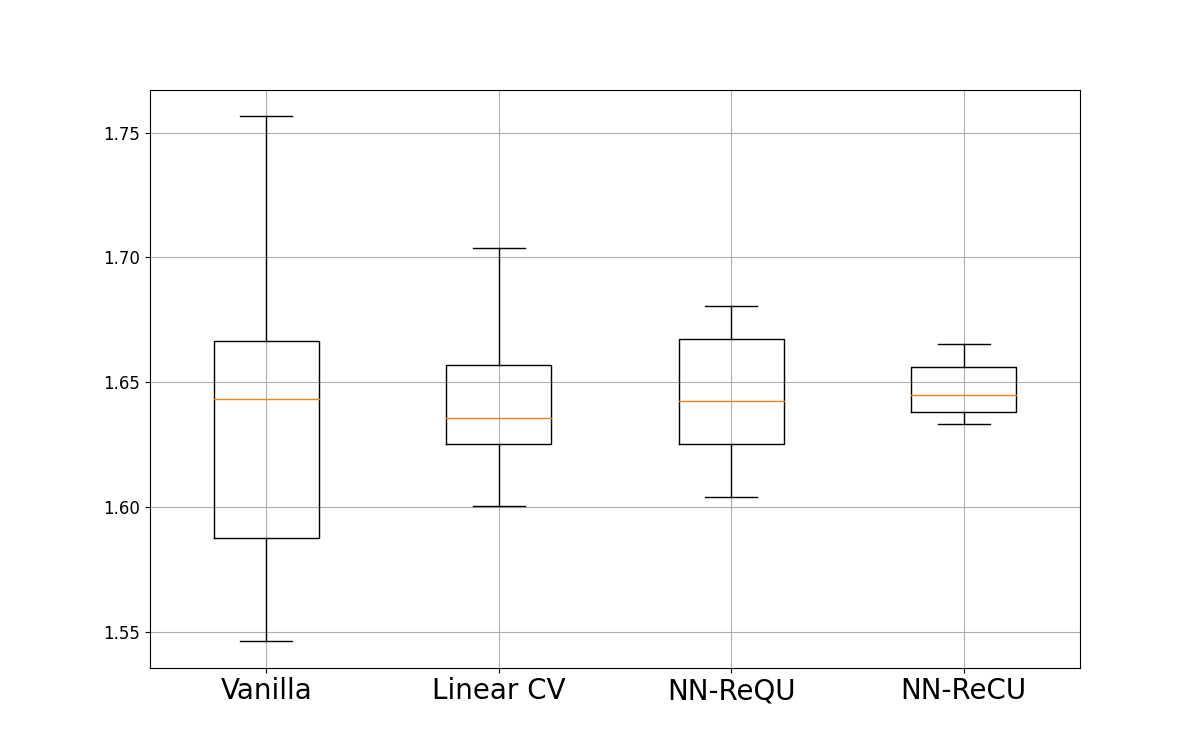}
      }
      \subfloat[Training chain]
      {
      \includegraphics[width=0.48\textwidth]{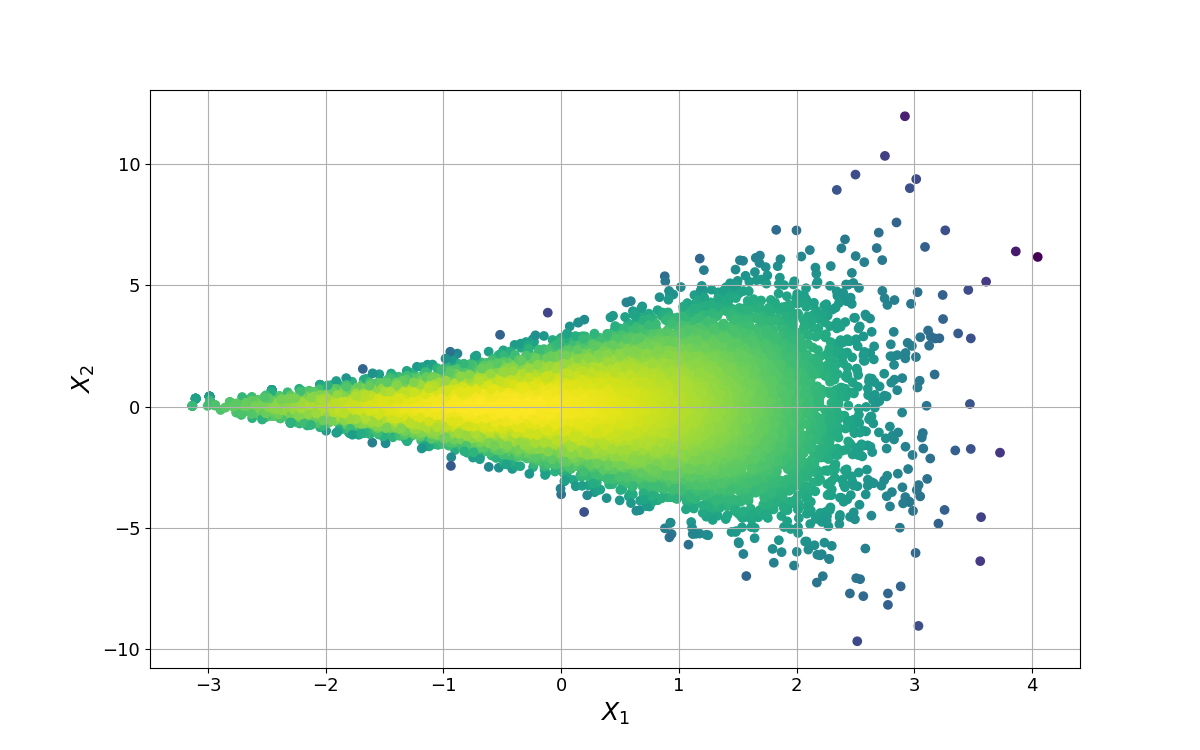}
      }
      
      \caption{Funnel distribution}
    \end{figure}
	
	\subsection{Banana-shaped distribution}
	
\begin{figure}[H]
      \centering
      \subfloat[Boxplots]{\includegraphics[width=0.48\textwidth]{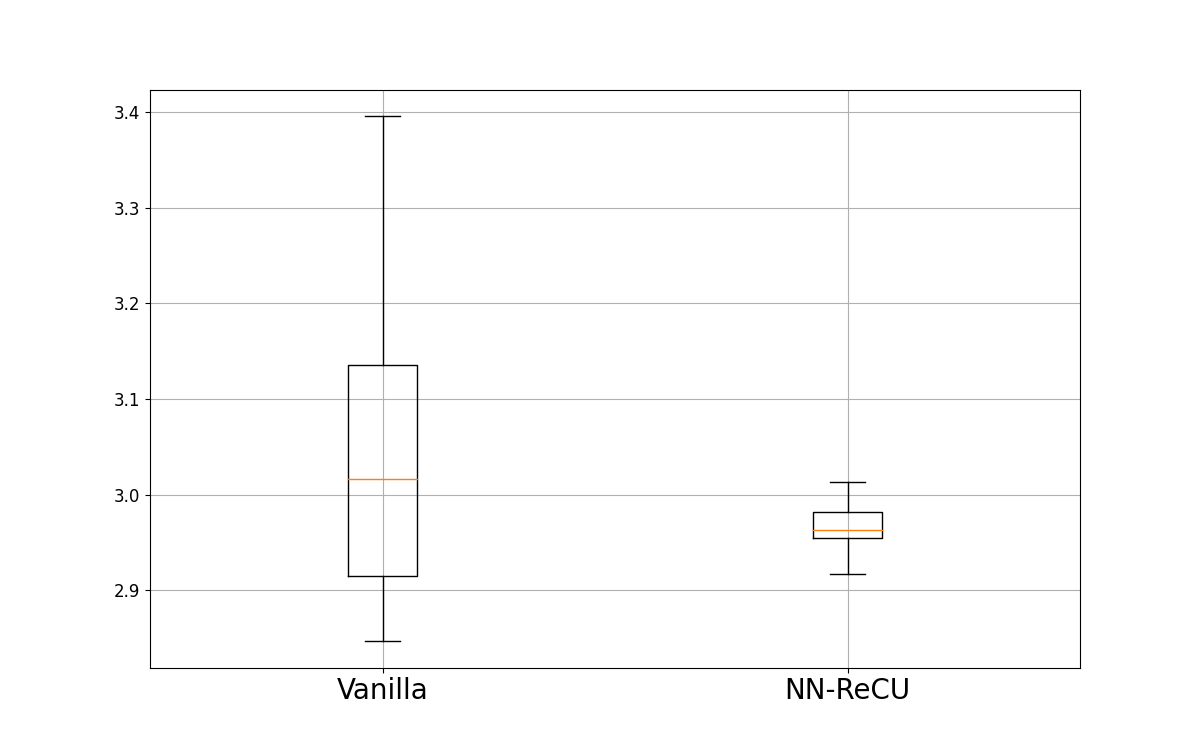}}
      \subfloat[Training chain, projection on first two coordinates]
      {
      \includegraphics[width=0.48\textwidth]{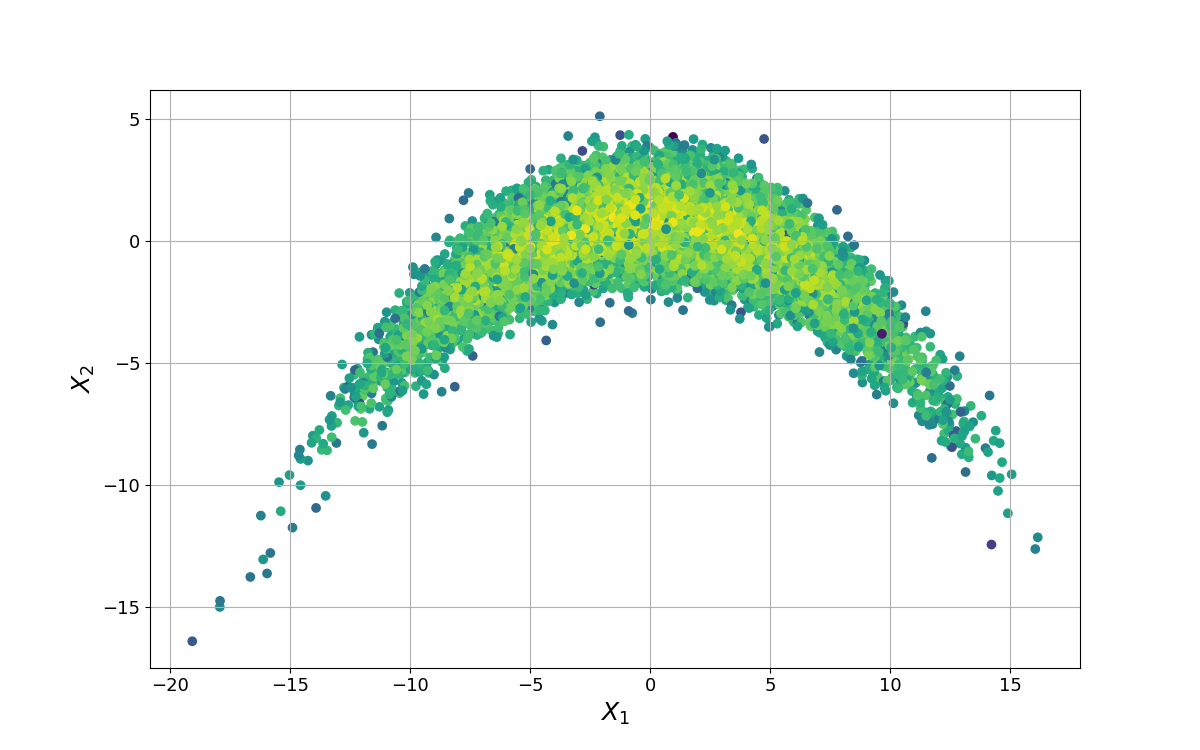}
      }
      \caption{Banana-shaped}
      \label{fig:banana}
    \end{figure}

	The ``Banana-shaped'' distribution is described in \cite{belomestny_variance_2020_esvm}. Its potential function is given by
	
	\begin{equation*}
	    U(x_1, \dots, x_d)
	    =
	    x_1^2/2p+(x_2+bx_1^2-pb)^2/2
	    +
	    \smm{k=3}{d}x_k^2/2
	\end{equation*}
	where parameters $p,b$ determine the curvature  and thickness of the density. We take $d=6,$ $p=20,$ $b=0.05.$
	The results  presented in Figure~\ref{fig:banana} and Table~\ref{tab:res-banana} are obtained using hyperparameters in Table~\ref{tab:hp-banana}. 
   \begin{table}[H]
    \centering
        \begin{tabular}{ |c|c|c|c|c|c|c| } 
        \hline
        CGM & 
        $\gamma$&
        $n_{\text{burn}}$& 
        $n_{\text{train}}$&
        $n_{\text{test}}$&
        $T$&
        $b_n$\\
        \hline
        ULA & 
        0.01&
        $10^5$& 
        $2\cdot 10^4$&
        $10^4$&
        30&
        30\\
        \hline
        \end{tabular}
        \caption{\label{tab:hp-banana}Banana-shaped, algorithm and chain generation hyperparameters}
    \end{table}

    \begin{table}[H]
    \centering
        \begin{tabular}{ |c|c| } 
        \hline
        Setup&
        $V_n(f)/V_n(f-g)$\\
        \hline
        ReCU&
        28\\
        \hline
        \end{tabular}
        \caption{\label{tab:res-banana} Banana-shaped, results}
    \end{table}
    
    \subsection{Logistic regression}

    We follow experimental setup from \cite{belomestny_variance_2020_esvm}.
    Let $\sY=(\sY_1,\ldots,\sY_n)\in \{0,1\}^n$ be a vector of binary response variables, $x \in\rset^d$ be a vector of regression coefficients, and $\sZ \in \rset^{N \times d}$ be a design matrix. 
    The log-likelihood and likelihood of $i$-th point for the logistic regression is given by
    \begin{align*}
      \logl{l}(\sY_i | x, \sZ_i) & = \sY_i  \sZ_i^{\sT} x -  \ln(1+\rme^{\sZ_i^{\sT} x}), \quad \likel{l}(\sY_i | x,\sZ_i) = \exp(\logl{l}(\sY_i | x, \sZ_i)), 
     \end{align*}
    where $\sZ_i^{\sT}$ is the $i$-th row of $\sZ$ for $i \in \lbrace1,\ldots,N\rbrace$. We complete the Bayesian model by considering the Zellner $g$-prior for the regression parameter $x$, that is, \(\mathcal{N}_{d}(0,g(\sZ^{\sT}\sZ)^{-1})\). Defining $\tilde{x}= (\sZ^{\sT} \sZ)^{1/2} x$ and $\tilde{\sZ}_i= (\sZ^{\sT} \sZ)^{-1/2} \sZ_i$, the scalar product is preserved, that is $\langle x, \sZ_i \rangle = \langle \tilde{x}, \tilde{\sZ}_i \rangle$ and, under the Zellner $g$-prior, $\tilde{x} \sim \mathcal{N}_{d}(0,g I_d)$. In the sequel, we apply the algorithms in the transformed parameter space with normalized covariates and put $g = 100$.
    \par
    The unnormalized posterior probability distribution $\pib{l}$ for the logistic regression model is defined for all $\tilde{x} \in \rset^d$ by
    \begin{align*}
     \pib{l}(\tilde{x} | \sY, \sZ) &\propto \exp(-\Ub{l}(\tilde{x})) \quad \text{with} \quad \Ub{l}(\tilde{x}) = -\sum\nolimits_{i=1}^{N}\logl{l}(\sY_i | \tilde{x}, \sZ_i) + (2g)^{-1}\norm[2]{\tilde{x}} .
     \end{align*}
    \par
    We analyze the performance of our algorithm on a Pima\footnote{\url{https://www.kaggle.com/uciml/pima-indians-diabetes-database}} dataset from the UCI repository. It contains  $N = 768$ observations in dimension $d = 9$ and their label value. We split the dataset into a training part $\mathcal{T}_{N-K}^{\text{train}}= [(y_i,\sZ_i)]_{i=1}^{N-K}$ and
    a test part $\mathcal{T}^{\text{test}}_K= [(y'_i,\sZ'_i)]_{i=1}^{K}$ by randomly picking $K$ test points from the data. In our experiment we pick $K=154$ (e.g. train part has size $\lfloor0.8\cdot N\rfloor=614$). Then we use ULA algorithm to sample from $\pib{l}(\tilde{x} | \sY, \sZ)$ and $ \pib{p}(\tilde{x} | \sY, \sZ)$ respectively. 
    \par
\begin{figure}[H]
      \centering
      \subfloat[Tanh, boxplots]
      {
      \includegraphics[width=0.48\textwidth]{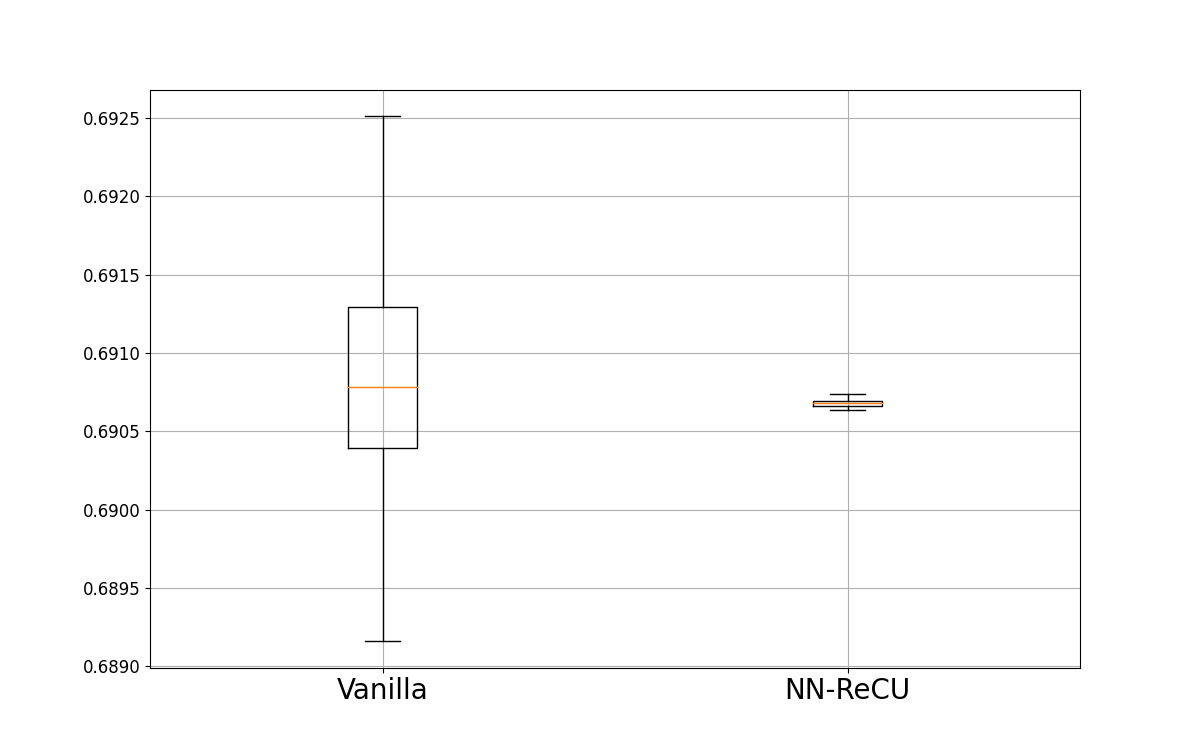}
      }
      \subfloat[Training chain, projection on first two coordinates]
      {
      \includegraphics[width=0.48\textwidth]{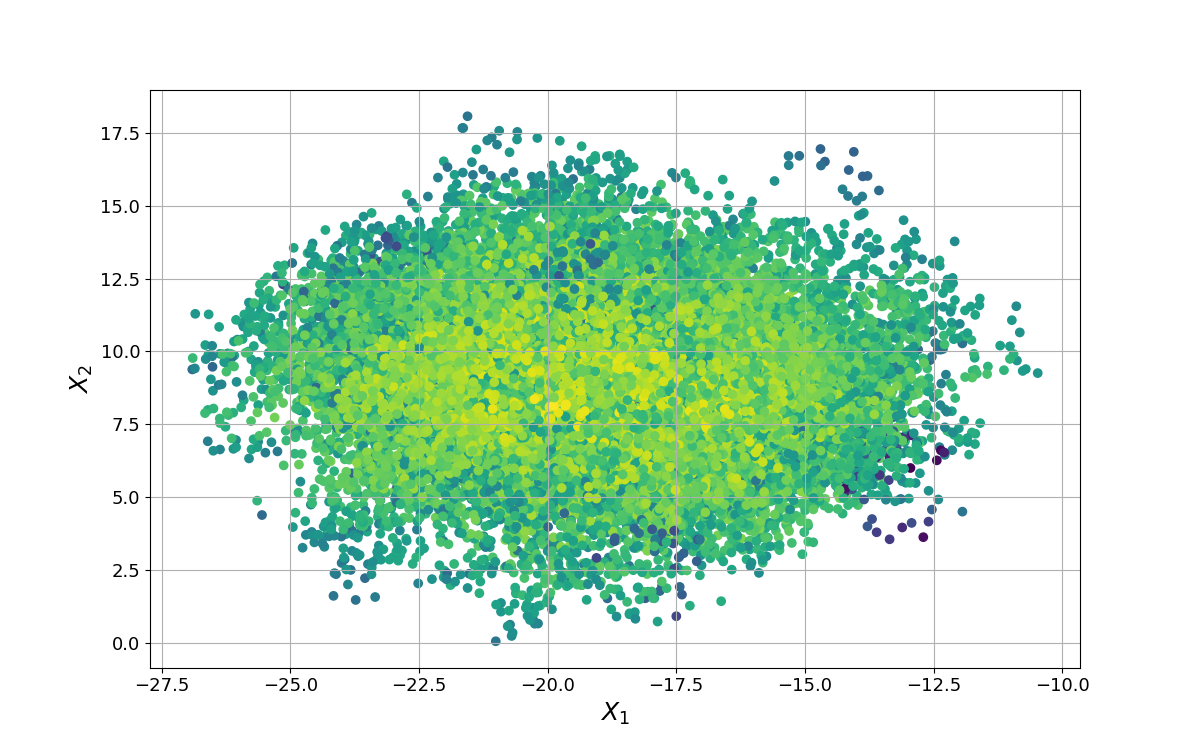}
      }
      
      \caption{\label{fig:lg-res}Logistic regression, Pima dataset}
    \end{figure}

    Given the sample $(\tilde{x}_k)_{k=0}^{n-1}$, we aim at estimating the average likelihood over the test set $\mathcal{T}^{\text{test}}_K$, that is,
    \[
    \int\nolimits_{\mathbb{R}^d}f(\tilde{x})\pib{l}(\tilde{x} | \sY, \sZ) \, \rmd \tilde{x} ,
    \]
    where the function $f$ is given by
    \[
    f(\tilde{x})= K^{-1} \sum_{i=1}^{K}\likel{l}(y'_i | \sZ'_i,\tilde{x}) .
    \]
    The results of ESVRR are reported in Table~\ref{tab:lr-res} and Figure~\ref{fig:lg-res} for 1-layered MLP with Tanh activation function.

    \begin{table}[H]
    \centering
        \begin{tabular}{ |c|c| } 
        \hline
        Setup&
        $V_n(f)/V_n(f-g)$\\
        \hline
        Tanh&
        122\\
        \hline
        \end{tabular}
        \caption{\label{tab:lr-res} Logistic Regression, results}
    \end{table}

    In this instance, it's noteworthy to observe a departure from our theoretical results by opting for the Tanh activation function instead of the ReCU activation function. Despite this deviation, the algorithm demonstrates effective performance and achieves notable variance reduction. This example underscores the flexibility in activation function choices, revealing that while we have established theoretical guarantees for ReCU-based neural networks, we are not strictly confined to them. There is potential to explore and establish bounds for the approximation error term associated with other activation functions, expanding the scope of applicability beyond the initially considered scenario.

	\section{Conclusion}
    \label{sec:conclusion}
    
    This research presents a novel theoretical solution to the challenge of variance reduction using a class of control variates based on neural networks. This methodology holds promise for diverse applications in machine learning, particularly within the realm of Bayesian statistics \cite{dellaportas2012control, mira2013zero}. A noteworthy aspect of this study is its reliance on contemporary findings regarding the approximation capabilities of neural networks \cite{schmidt-hieber2020, belomestny2022simultaneous}. The conducted numerical experiments highlight a substantial reduction in variance when applying the proposed algorithm to various datasets.  These experiments also suggest the importance of using neural networks with smooth activation functions for Stein control functionals. There are several avenues for extending this research. A comprehensive numerical investigation, inspired by \cite{si_scalable_2020}, could delve into implementation intricacies, potentially providing deeper insights into the relationship between optimization, variance convergence, and associated properties, while considering trade-offs in terms of time, memory, and quality. Furthermore, future research could explore broader function classes beyond those constrained by the current assumptions. Lastly, the outlined approach might be adapted for variance reduction methods beyond Stein control variates or for applications in solving partial differential equations using neural networks.
    
    \section{Acknowledgements}
    This project was supported by the grant for research centers in the field of AI provided by the Analytical Center for the Government of the Russian Federation (ACRF) in accordance with the agreement on the provision of subsidies (identifier of the agreement 000000D730321P5Q0002) and the agreement with HSE University No. 70-2021-00139.

    \section{Proofs}
     \label{sec: proofres}
    \subsection{Proof of Theorem~\ref{th:bound_approx}}
    \paragraph*{\textbf{Step 1}} We aim to trace the precise dependence of our bounds upon the factors $n,R$, and $K$. Multiplicative dependencies upon the remaining model parameters are accounted in $\lesssim$. From now on we assume that the parameters $R,K>0$ are fixed, and later we optimize them over the number of samples $n$. Recall that the underlying class of control variates $\Phi_{R}$ is given in \Cref{th:nn_constr}. For notation simplicity we denote
    \[
    \Phi\eqdef \Phi_{R} ,\quad \H \eqdef
    \mathcal{H}_{\Phi_R}\,,
    \]
    omitting the explicit dependence of both classes upon $R$. Now we bound $\AV{\hhf{\varphi}{R}}$ for $\varphi\in \Phi$. Note that the functions in $\H$ are bounded, due to the construction of $\Phi$ from \eqref{eq:phi_def}. Let us denote by $B_{\H}$ the upper bound for $\sup_{h\in\H}{|h|_{\infty}}$, that is, $\sup_{h\in\H}{|h|_{\infty}} \leq B_{\H}$. Thus we get using \Cref{th:main_slow} that for each $x_0 \in \Xset_0$ 
 and $\delta \in (0;1)$ there exists $n_0=n_0(\delta)>0$ such that for all $n\geq n_0$ it holds with $\P_{x_0}-$probability at least $1-\delta$ that
    \begin{multline}
    \label{eq:bound_mcmc}
        \AV{\hh{n, \eps,R}}
        =
        \AV{\hh{n, \eps,R}}
        -
        \inf_{\varphi\in{\Phi}}\AV{\hhf{\varphi}{R} }
        +
        \inf_{\varphi\in{\Phi}}\AV{\hhf{\varphi}{R}}\\
        \lesssim
        \inf_{\varphi\in{\Phi}}\AV{\hhf{\varphi}{R}}
        +
        B_{\H}^2
        \log(n) \gamma_{\ltwo(\pi)}(\H,n)
		+ (B_{\H}\vee \sup_{h\in\H}\norm{h}_{W^{1/2}}^2)^2\frac{\log(n)\log(\fraca{1}{\delta})}{\sqrt{n}}.
    \end{multline}
    In the bound above we used that $\inf_{h\in{\H}}\AV{h} = \inf_{\varphi\in{\Phi}}\AV{\hhf{\varphi}{R}}$. Now we bound $\inf_{\varphi\in{\Phi}}\AV{\hhf{\varphi}{R}}$. Obviously, infimum is upper-bounded by value of $V_\infty$ on any $\hhf{\varphi}{R}\in\H$. We pick $\hhf{\psi}{R}$ which corresponds to $\psi \in \Phi$, such that $\norm{\psi - \varphi^*_R}_{\H^3(\Aset)}\leq \epsnnapprox$. The existence of such $\psi$ is guaranteed by \Cref{th:nn_constr}. Thus, we get from \Cref{th:var_decomp} that
    \begin{align*}
    \inf_{\varphi\in{\Phi}}\AV{\hhf{\varphi}{R}} \lesssim r \norm{\thhf{\psi}{R}}_{L^2(\Xset, \pi)}
        \bigl(
        (1+\pg(R)) d \epsnnapprox + \norm{\tilde{f}}_{L^2(\setcomp{\Aset}, \pi)}
        \bigr) +\norm{\thhf{\psi}{R}}^2_{W^{1/2}} \frac{\rho^{(r+1)/2}}{1-\rho^{1/2}}\eqsp.
	\end{align*}
    Here we additionally used that $\| \cdot \|_{C^2(A)}\leq \| \cdot \|_{\H^3(A)}$. 
    \paragraph*{\textbf{Step 2}} Now we need to quantify the dependence of covering number of the class $\H$ on $K$ and $H_{R}$ from \Cref{prop:recu_approx_const}. To do this, we note that covering number for $\H = \mathcal{H}_{\Phi_R, \eps} $ is not greater than covering number for $\mathcal{H}_{\NN}$, where $\NN=\NN(L, \A(K, H_{R}), \snzw(K, H_{R}))$. Moreover, for $f, h \in \H$, such that $\|f-h\|_{L^\infty(\Xset)}\leq \eps$, it is easy to check that 
    $\|f-h\|_{L^2(\Xset, \pi)} \leq \|f- h\|_{L^\infty(\Xset)}\leq \eps$.
    Thus, an $\eps$-covering net in $L^\infty (\Xset)$ is also an $\eps$-covering net for $L^2(\Xset, \pi)$, and $\gamma_{\ltwo(\pi)}(\mathcal{H},n) = \gamma_{\ltwo(\Xset, \pi)}(\mathcal{H},n) \leq \gamma_{L^\infty(\Xset)}(\mathcal{H},n)$. 
    \par 
    Thus, with the definition \eqref{eq:f_minus_steincv} of $\H_\NN$, it is enough to bound the covering number of $\H_\NN$ in $C^2$ metric, which can be controlled with \Cref{Lem:h2_norm_bound}. Indeed,
    \[
    \mathcal{N}(\NN, \|\cdot\|_{C^2([0,1]^{d})}, \eps')
    \leq \left(2(\eps')^{-1} L^2(L+1)(L+2)3^{3L+4}
V^{L3^{L+5}+1}\right)^{\snzw+1}\,,
    \]
    where $V$ is defined in \eqref{eq:V_def}, $L$ is the depth of $\NN$, and $q$ is the number of non-zero weights. Using \Cref{prop:recu_approx_const}, we get
    \begin{align*}
        L \leq 8\cdot(4+2(\lbeta-2)+\lceil\log_2 d\rceil + 3),\;
        q\lesssim K^dR^{\hpow+3} H_{R},\;
        V\lesssim (K^dR^{\hpow+3} H_{R})^L\,.
    \end{align*}
    Moreover, by \Cref{lem:cv_diff_bound}, 
    \[
    \|f-h\|_{L^2([0,1]^{p_0})}
    \leq 
    (1+\pg(R))d\|f-h\|_{C^2([0,1]^{p_0})}\,,
    \]
    hence
    \begin{align*}
    \log\mathcal{N}(\mathcal{H}_\NN, \|\cdot\|_{L^2([0,1]^{d})}, \eps) 
    &\leq \log \mathcal{N}(\NN, \|\cdot\|_{C^2([0,1]^{d})}, \frac{\eps}{(1+\pg(R))d}) \\ 
    &\leq (\snzw+1)\log(2(1+\pg(R))d\eps^{-1} L^2(L+1)(L+2)3^{3L+4}
V^{L3^{L+5}+1}) \\
    &\lesssim K^dR^{\hpow+3} H_{R}
    \left(
    \log((1+\pg(R))\eps^{-1})
    +
    \log (K^dR^{\hpow+3} H_{R})
    \right).
    \end{align*}
    Thus, $\gamma_{\ltwo(\pi)}(\mathcal{H},n)$ is the smallest number $\eps$ satisfying
    \begin{equation}
    \label{eq:cov_num_ineq}
    K^dR^{\hpow+3} H_{R}
    \left(
    \log((1+\pg(R))\eps^{-1})
    +
    \log (K^dR^{\hpow+3} H_{R})
    \right)
    \leq n\eps^2.
    \end{equation}
    
    \paragraph*{\textbf{Step 3}}
    Now we substitute constructed class with chosen approximator into the bound:    
    \begin{multline*}
		\AV{\hh{n, \eps,R}} 
		\lesssim
		r  \norm{\thhf{\psi}{R}}_{L^2(\Xset, \pi)}
    	(
        (1+\pg(R)) d \left(
		\frac{R^{\hpow+3}H_{R}}{K^{\hpow}}
		\right)
        +
        \|\tilde{f}\|_{L^2(\setcomp{\Aset}, \pi)}
        )\\
         + \|\thhf{\psi}{R}\|^2_{W^{1/2}}
         \frac{\rho^{(r+1)/2}}{1-\rho^{1/2}} 
		+
		(B_{\H}\vee \sup_{h\in\H}\norm{h}_{W^{1/2}}^2)^2\ 
        \biggl( \log(n) \gamma_{\ltwo(\pi)}(\mathcal{H},n)
		+ \frac{\log(n)\log(\fraca{1}{\delta})}{\sqrt{n}}\biggr).
	\end{multline*}	
	Firstly, from the fixed point inequality \eqref{eq:cov_num_ineq} it is clear, that 
	\[
	\gamma_{\ltwo(\pi)}(\mathcal{H},n)=\eps \leq \sqrt{\frac{K^dR^{\hpow+3} H_{R}}{n}}\log n
	\]
    Now recall, that $\psi \in \Phi$ satisfies \eqref{eq:close_net_prop} for given $K, R$. By \Cref{th:nn_constr}, we have $\psi \in \mathcal{H}^{3}(\Aset, (2R)^{\hpow+3}H_{R} + \epsnnapprox)$, where $\epsnnapprox \lesssim \frac{C^{(\hpow+3) d} H_{R} (\hpow+3)^3}{K^{\hpow}}$. Since the function $\psi$ is bounded, 
	\[
	\norm{\thhf{\psi}{R}}_{L^2(\Xset, \pi)}
	\lesssim B_{\H},\;
	\norm{\thhf{\psi}{R}}_{W^{1/2}}^2 \lesssim B_{\H}^2, \quad \sup_{h\in\H}\norm{h}_{W^{1/2}} \lesssim B_{\H}\,.
	\]
    Now we pick $B_{\H}$ which bounds every $f-\ggf{\psi}{R}=\hhf{\psi}{R}\in\H$. Indeed,
    \[
    \norm{f-\ggf{\psi}{R}}_{L^\infty(\Xset)}
            \leq
            \norm{f}_{L^\infty(\Xset)}
            +
            \norm{\ggf{\psi}{R}}_{L^\infty(\Xset)}
            \lesssim B_f
            +
            (1+\pg(R))d \|\psi\|_{C^2(\Aset)}\,,
    \]
    where $B_f$ is from \Cref{assu:AUF}, and the latter inequality is due to \Cref{lem:cv_diff_bound}. Since $\psi \in \Phi$,
    \[
    \|\psi\|_{C^2(\Aset)} \leq \|\psi\|_{\H^3(\Aset)} \leq 2((2R)^{\hpow+3}H_{R} + \epsnnapprox)
    \lesssim 
    (2R)^{\hpow+3}H_{R} + \frac{ R^{(\hpow+3)} H_{R} (\hpow+3)^3}{K^{\hpow}}.
    \]
    Next, using \Cref{lem:gamma_func_beh}, we get
    \begin{equation}
    \label{eq:outer_bound}
    \|\tilde{f}\|_{L^2(\setcomp{\Aset}, \pi)} \lesssim R^d e^{-R^2}\,.
    \end{equation}
    Recall also, that  $\pg(R)\leq C_0\cdot R^{\udeg}$ for $R\geq R_0$ by \Cref{assu:AUF}. Now we need to account for the dependence of the H\"older constant $H_{R}$ of $\varphi^*_R$ upon $R$. This dependence is subject to the H\"older constants of elliptic operator coefficients in equation \eqref{eq:dirichletarray}, which change with growth of $R$, and to the growing domain $\Aset$. Applying \Cref{lem:sol_holder_const}, we get that $H_{R} \lesssim R^{\hdeg}$, where the power $\hdeg$ defined in \Cref{lem:sol_holder_const} is independent of the problem dimension $d$. Combining the bounds above, 
\begin{multline}\label{eq:thmbound3}
\AV{\hh{n, \eps,R}} \lesssim R^{5(\hpow+3+\udeg)+\hdeg}
        \biggl(
		r \biggl(
		\frac{R^{\hpow+3+\hdeg}}{K^{\hpow}} + R^de^{-R^2}
        \biggr)
         + 
         \frac{\rho^{(r+1)/2}}{1-\rho^{1/2}} \\
		+\log^2(n) \sqrt{\frac{K^{d}R^{\hdeg+\hpow+3}}{n}}
		+ \frac{\log(n)\log(\fraca{1}{\delta})}{\sqrt{n}}
		\biggr).
\end{multline}	
Now it remains to optimize the right-hand side of \eqref{eq:thmbound3} over $n$, choosing the values $r$, $R$, and $K$. We select $K=n^\alpha$ for some $\alpha > 0$, $R=\log n$, $r = -\log_\rho n$ ($0<\rho<1$), and note that the bound \eqref{eq:cov_num_ineq} holds for $n \geq n_0$. Optimising over the power $\alpha$ in \eqref{eq:thmbound3}, we end up with the choice 	
\[
\textstyle{K = n^{\frac{1}{2\hpow+d}}\eqsp, \quad \eps = \sqrt{K^d/n} = n^{-\frac{\hpow}{2\hpow+d}}}\eqsp,
\]
and the final bound yields
\[
\AV{\hh{n, \eps,R}} \lesssim (\log n)^{C_{\ref{th:bound_approx}, 1}} \bigl( n^{-\frac{\hpow}{2\hpow+d}} \vee \frac{\log(1/\delta)}{\sqrt{n}} \bigr)\eqsp,
\]
	where $C_{\ref{th:bound_approx}, 1}=C_{\ref{th:bound_approx}, 1}(\hpow, \udeg, \hdeg)$. Note, that the constant $C_{\ref{th:bound_approx}, 1}$ is independent of dimensionality $d$. 

\paragraph*{\textbf{Step 4}} Denote by $\hat{g}_{n,\eps,R}$ the truncated control variate corresponding to the function $\hh{n, \eps,R}$, that is, $\hat{g}_{n,\eps,R} = f - \hh{n, \eps,R}$. Note that in general $\pi(\hat{g}_{n,\eps,R}) \neq 0$ due to the truncation. Let us see how this bias can be accounted in the final bound for $\pi(f)$. The CLT yields the asymptotic $1-\delta$ confidence interval for $\pi(\hh{n, \eps,R})$ of the form
\[
        \pi(\hh{n, \eps,R})\in
    	\biggl(
    	    \pi_N(\hh{n, \eps,R}) - 
    	    c_{1-\delta/2}N^{-1/2}\sqrt{\AV{\hh{n, \eps,R}}};
    	    \pi_N(\hh{n, \eps,R}) +c_{1-\delta/2}  N^{-1/2}\sqrt{\AV{\hh{n, \eps,R}}}
    	\biggr)\eqsp.
\]
Recall, that by \eqref{eq:ggf} $\ggf{\varphi}{R} = g_\varphi \indibr{x\in \Aset}$, where $\pi(g_\varphi)=0$. Thus,    
\begin{equation}
\label{eq:bias_g_phi}
0 = \pi(g_\varphi) = \pi(\ggf{\varphi}{R})+\pi(g_{\varphi} \indibr{\setcomp{\Aset}})\eqsp.
\end{equation}
Denote by $\varphi_{n,\eps}$ the function, which corresponds to $\hat{g}_{n,\eps,R}$ in \eqref{eq:ggf}, that is, 
\[
\hat{g}_{n,\eps,R}(x) = \{\Delta \varphi_{n,\eps}(x) +
\langle \nabla\log \pi(x), \nabla \varphi_{n,\eps}(x) \rangle\} \indibr{x\in \Aset}\eqsp.
\]
Then, with \eqref{eq:bias_g_phi}, we get an asymptotic $1-\delta$ interval for $\pi(f)$ of the form  
\begin{align*}
\pi(f) \in
\bigl[ 
\pi_N(\hh{n, \eps,R}) -\pi(g_{\varphi_{n,\eps}} \indibr{\setcomp{\Aset}})  - c_{1-\delta/2}N^{-1/2}\sqrt{\AV{\hh{n, \eps,R}}}; \\
\pi_N(\hh{n, \eps,R}) -\pi(g_{\varphi_{n,\eps}} \indibr{\setcomp{\Aset}}) + c_{1-\delta/2}N^{-1/2}\sqrt{\AV{\hh{n, \eps,R}}}
\bigr]\eqsp.
\end{align*}
Now we bound the appearing bias term $\pi(g_{\varphi_{n,\eps}} \indibr{\setcomp{\Aset}})$. It is clear that the bias decreases with the increase of the truncation radius $R$, which explains why we have to choose $R$ increasing with the sample size $n$. Recall that $\varphi_{n,\eps} \in \Phi$, where the class $\Phi$ is constructed in \Cref{th:nn_constr}. Combining  \Cref{lem:cv_diff_bound} and \Cref{lem:infty_norm_hessian}, we get 
    
\begin{align*}
\norm{g_{\varphi}}_{L^\infty(\Aset)} & \lesssim 
    (1+\pg(R))d \|\varphi\|_{C^2(\Aset)}
    \lesssim
    (1+\pg(R))d
    \bigl\{\prod_{\ell=0}^{L}(p_{\ell}+1)^{4(L3^{L}+1)+1}\bigr\}\\
    &\lesssim 
    R^\udeg (K^dR^{\hpow+3} H)^{L^23^{L+2}} 
    \lesssim 
    R^{\udeg + L^23^{L+2}(\hpow+3+\hdeg)} K^{dL^23^{L+2}}.
\end{align*}
    Note that due to rescaling from $\Aset$ to $[0;1]^d$, neural network has the latter as its domain.
    Since $g_\varphi(x)$ grows not faster than a polynomial in $R$ for $\|x\| = R$, we get from \Cref{lem:gamma_func_beh} that
    \begin{equation}
    \label{eq:bias_bound}
|\pi(g_{\varphi}\cdot\indibr{\setcomp{\Aset}})|\leq 
    \|g_{\varphi}\cdot\indibr{\setcomp{\Aset}}\|^2_{L^2(\setcomp{\Aset}, \pi)}
    \lesssim 
    R^{d+2(\udeg + L^23^{L+2}(\hpow+3+\hdeg))} K^{dL^23^{L+2}}
    e^{-R^2}.
    \end{equation}
    Thus, for the choice $R = \log (n)$, $K = n^{\frac{1}{2\hpow+d}}$, we get 
    \[
    |\pi(g_{\varphi}\cdot\indibr{\setcomp{\Aset}})|
    \lesssim 
    (\log n)^{C_{\ref{th:bound_approx}, 2}}\frac{n^{C_{\ref{th:bound_approx}, 3}}}{n^{\log n}},
    \]
    where $C_{\ref{th:bound_approx}, 2} = C_{\ref{th:bound_approx}, 2}(\hpow, d, \udeg)$, and $C_{\ref{th:bound_approx}, 3} = C_{\ref{th:bound_approx}, 3}(\hpow, d)$. Combining the bounds above yields the statement of the theorem.

    \subsection{Auxiliary results}

    \begin{Lem}
    \label{th:var_decomp}
    Assume \Cref{assu:AUF}, \Cref{assu:ge}, and let $\varphi\in\H^{\hpow+3}(\Aset, H_{R})$ for $\hpow>0$ and $R>0$. Then for any $r\in\nset$ it holds that
    \begin{multline}
    \label{eq:th_approx_err_bound}
        \AV{\hhf{\varphi}{R}} \lesssim 
        r
        \norm{\thhf{\varphi}{R}}_{L^2(\Xset, \pi)}\biggl(
        (1+\pg(R)) d \|\varphi-\varphi^*_R\|_{C^2(\Aset)}
        +
        \norm{\tilde{f}}_{L^2(\setcomp{\Aset}, \pi)}
        \biggr)\\
        +
        \norm{\thhf{\varphi}{R}}^2_{W^{1/2}}
        \frac{\rho^{(r+1)/2}}{1-\rho^{1/2}}\eqsp.
    \end{multline} 
    Here $\lesssim$ stands for an inequality up to a constant $C(\varsigma, \pi(W))$, $\pg$ is defined in \Cref{assu:AUF}, and $\varphi^*_R$ is from \ref{eq:diff_eq_alt}.
    \end{Lem}
    \begin{proof}
    Applying \Cref{lem:cov_sum_bound} to the functions $\hhf{\varphi}{R}$ and $\hhf{\varphi_R^*}{R}$, we get
    \begin{multline*}
        \AV{\hhf{\varphi}{R}} \leq (2r+1) \sqrt{\PVar{\pi}(\hhf{\varphi}{R})} \left(
        \sqrt{\PVar{\pi}(\hhf{\varphi}{R}-\hhf{\varphi^*_R}{R})}
        +
        \sqrt{\PVar{\pi}(\hhf{\varphi^*_R}{R})}
        \right)\\
         + 2
        \left(
        \varsigma^{1/2}\pi(W)\|\thhf{\varphi}{R}\|^2_{W^{1/2}}
        \right)
        \frac{\rho^{(r+1)/2}}{1-\rho^{1/2}}\eqsp.
    \end{multline*} 
Bounding the variance term $\PVar{\pi}(\hhf{\varphi}{R}-\hhf{\varphi^*_R}{R})$, we get
    \begin{align*}
    \sqrt{\PVar{\pi}(\hhf{\varphi}{R}-\hhf{\varphi^*_R}{R})}
    \leq
    \norm{\hhf{\varphi}{R}-\hhf{\varphi^*_R}{R}}_{L^2(\Xset, \pi)}=
    \norm{\ggf{\varphi^*_R}{R}-\ggf{\varphi}{R}}_{L^2(\Aset, \pi)}\,,
    \end{align*}
    where the first inequality is due to an obvious bound $\Var[X]\leq \E[X^2]$.
    \par
    Now we aim to apply \Cref{lem:cv_diff_bound} to control $\norm{\ggf{\varphi^*_R}{R}-\ggf{\varphi}{R}}_{L^2(\Aset, \pi)}$. It requires to ensure that $\varphi^*_R \in \mathcal{H}^{\hpow+3}(\Aset, H_{R}),\ \hpow>0$, for some $H_{R} > 0$.  Recall that $\varphi^*_R$ is a solution of the PDE \eqref{eq:diff_eq_alt}. Thus, under \Cref{assu:AUF} we apply \Cref{lem:sol_holder_const}, and obtain that $\varphi^*_R \in \mathcal{H}^{\hpow+3}(\Aset, H_{R})$ for some constant $H_{R} > 0$, which possibly depends upon $R$. Thus, the application of \Cref{lem:cv_diff_bound} yields
    \[
    \norm{\ggf{\varphi^*_R}{R}-\ggf{\varphi}{R}}_{L^2(\Aset, \pi)} \leq (1+\pg(R))d \|\varphi - \varphi^*\|_{C^2(\Aset)}\,.
    \]
    Also note, that 
    \[
    \textstyle{
    \PVar{\pi}(\hhf{\varphi^*_R}{R})
        \leq
        \E_\pi
        [(\hhf{\varphi^*_R}{R}-\pi(f))^2]
        =
        \int_{\setcomp{\Aset}}
        (f(x)-\pi(f))^2 \pi(dx) = 
        \|\tilde{f}\|^{2}_{L^2(\setcomp{\Aset}, \pi)}}\eqsp.
    \]
    Combining the bounds above yields \eqref{eq:th_approx_err_bound}.
    \end{proof}

    \begin{Lem}
    \label{th:nn_constr} 
    Assume \Cref{assu:AUF}. Then for each $R, K > 0$, the class $\Phi_R$ defined in \eqref{eq:phi_def} is a non-empty class of twice-differentiable functions. Moreover, there exists $\psi \in \Phi_{R}$, such that
    \begin{equation}\label{eq:close_net_prop}
    \norm{\psi - \varphi^*_R}_{\H^3(\Aset)}\leq \epsnnapprox\,.
    \end{equation} 
     where $\epsnnapprox$ is defined in  \eqref{eq:eps_nnaprox_def}, and $\varphi^*_R$ is defined in \eqref{eq:diff_eq_alt}.
    \end{Lem}
    \begin{proof} Recall that under \Cref{assu:AUF}, the function $\varphi^*_R\in \mathcal{H}^{\hpow+3}(\Aset, H)$, $\hpow>0$. Thus we aim to approximate $\varphi^*_R$ within the ReCU neural networks, using the construction of \Cref{prop:recu_approx_const}. Thus it is enough to rescale the function $\varphi^*_R$ in order to change its domain from $[-R;R]^{d}$ to $[0,1]^{d}$. With a shift vector $v = \frac{1}{2}(1,1,\dots,1)^T$, the latter can be achieved with a transformation 
    \[
    \|f\bigl(t/(2R)+v\bigr) - \varphi^*_R(t)\|_{\H^3([-R;R]^d)} =
    \|f\left(x\right)-\varphi^*_R(2R(x-v))\|_{\H^3([0;1]^d)}.
    \]
    Thus, it is enough to approximate the function $\varphi^*_R(2R(x-v))$, $x \in [0,1]^{d}$. It is easy to check that 
    \[
    \varphi^*_R(2R(x-v))\in\H^{\hpow+3}([0;1]^d, (2R)^{\hpow+3}H).
    \]
    \Cref{prop:recu_approx_const} describes the function class $\NN(L, \A(K, H), \snzw(K, H))$ containing function $u$, such that 
    \begin{multline}
    \label{eq:nn_approx_property}
    \|u\bigl(x/(2R) + v\bigr)-\varphi^*_R(x)\|_{\H^3([-R;R]^d)}\leq 
    \|u(x)-\varphi^*_R(2R(x-v))\|_{\H^3([0;1]^d)}
    \\
    \leq \frac{C^{d(\hpow+3)}(2R)^{(\hpow+3)}H(\hpow+3)^3 }{K^{\hpow}} = \epsnnapprox.
    \end{multline}
    Note that functions from class $\NN(L,\A(H,K),\snzw(H, K))$ are twice continuously differentiable, because we use ReCU activation in neural networks. Due to \Cref{prop:recu_approx_const}, the constructed rescaled neural network $u\left(x /(2R) + v\right)$ lies in $\mathcal{H}^{s+3}(\Aset, (2R)^{\hpow+3}H + \epsnnapprox)$, where $s$ is the largest integer strictly less than $\hpow$, and $\epsnnapprox$ is from \eqref{eq:nn_approx_property}. 
    \par 
    Finally, we impose additional constraints on the class $\NN(L,\A(H,K),\snzw(H, K))$ and consider the functions which are close to the constructed mapping $u(1/(2R) x + v)$. Thus we consider class
    \begin{align*}
    \Phi_R
    =
    \biggl\{
    \varphi = g\circ
    \left(
    \frac{1}{2R}x+v
    \right)
    \biggm|
    \norm{
    \varphi
    }_{\H^3(\Aset)}\leq 2((2R)^{\hpow+3} H + \epsnnapprox
    )\,, \quad g \in \NN (L, \A(H, K),\snzw(H, K))
    \biggr\}\,.
    \end{align*}
    Note that at least the approximating function $u\left(x/(2R) + v\right)$ lies in this class, thus it is non-empty. Here the constant $2$ in the norm bound $2 ((2R)^{\hpow+3} H + \epsnnapprox)$ is picked only for convenience and can be substituted by $1 + \delta$ for any $\delta > 0$.     
    \end{proof}
    
Now we provide some simple lemmas required to bound the covariance under arbitrary initial distribution.

\begin{Lem}[\protect{\cite[Lemma~12]{belomestny_variance_2020_esvm}}]
\label{lem:cov_bound}
Assume \Cref{assu:ge}. Then for any $h\in \mathcal{H}$, $x \in \Xset$, and $s\in\nzeroset$,
\begin{equation*}
|E_x[\tilde{h}(X_0)\tilde{h}(X_s)]|
\leq \varsigma^{1/2}\rho^{s/2}W(x)\|\tilde{h}\|^2_{W^{1/2}},
\end{equation*}
and
\begin{equation*}
\label{eq:cov_bound}
|\rho_\pi^{(h)}(s)|
\leq \varsigma^{1/2}\rho^{s/2}\pi(W)\|\tilde{h}\|^2_{W^{1/2}}.
\end{equation*}
\end{Lem}

\begin{Lem}
\label{lem:cov_sum_bound}
Assume \Cref{assu:ge} and let 
$g, h \in L^2(\pi)$. Then, for any $r\in\nset$ it holds that
\begin{equation}
\label{eq:v_infty_bound}
V_\infty (g) \leq (2r+1) \sqrt{\PVar{\pi}(g)} \left(
\sqrt{\PVar{\pi}(g-h)}
+
\sqrt{\PVar{\pi}(h)}
\right) +
\frac{2 \varsigma^{1/2}\pi(W)\|\tilde{g}\|^2_{W^{1/2}} \rho^{(r+1)/2}}{1-\rho^{1/2}}\eqsp.
\end{equation}
\end{Lem}
\begin{proof}
Using the definition of the asymptotic variance $V_\infty(\cdot)$ in \eqref{eq:var_def_mcmc}, 
\begin{align*}
V_\infty (g) \leq \PVar{\pi} (g) + 2\smm{k=1}{\infty} \left|\PCov(g(X_0), g(X_k))\right|\eqsp.
\end{align*}
Now we split the sum above into $2$ parts. For $1 \leq k \leq r$, the covariances are bounded as
\begin{align*}
\bigl|\PCov(g(X_0), g(X_k))\bigr| 
&\leq \left|\PCov(g(X_0), g(X_k)-h(X_k))\right| +
\left|\PCov(g(X_0), h(X_k))\right| \\
&\leq
\sqrt{\PVar{\pi}(g)} \sqrt{\PVar{\pi}(g-h)}+ 
\sqrt{\PVar{\pi}(g)} \sqrt{\PVar{\pi}(h)}\,.
\end{align*}
To bound the remainder sum, we use \Cref{lem:cov_bound}, and obtain
\begin{align*}
2\smm{k=r+1}{\infty}
\left|\PCov(g(X_0), g(X_k))\right| \leq 2 \varsigma^{1/2}\pi(W)\|\tilde{g}\|^2_{W^{1/2}} \smm{k=r+1}{\infty}\rho^{k/2}  =
2\varsigma^{1/2}\pi(W)\|\tilde{g}\|^2_{W^{1/2}}\frac{\rho^{(r+1)/2}}{1-\rho^{1/2}}\eqsp
\end{align*}
Combining the bounds above yields \eqref{eq:v_infty_bound}.
\end{proof}

Next we aim to replace $\|.\|_{L^2(\Aset, \pi)}$ with $\|.\|_{H^{3}(\Aset)}$ on $\Aset$. 
\begin{Lem}
\label{lem:cv_diff_bound}
Assume \Cref{assu:AUF} and $\varphi,\ \varphi^* \in \mathcal{H}^{\hpow+3}(\Aset, H)$, $\hpow>0$. Then 
\[
\|\ggf{\varphi -\varphi^*}{R}\|_{L^2(\Aset, \pi)}=
\|\ggf{\varphi}{R} -\ggf{\varphi^*}{R}\|_{L^2(\Aset, \pi)} \leq (1+\pg(R))d \|\varphi - \varphi^*\|_{C^2(\Aset)}.
\]
\end{Lem}
\begin{proof}
Recall that $\pi$ is a probability measure on $\Xset$. Thus, with the definition \eqref{eq:ggf} and H\"older's inequality, 
\begin{align*}
\|\ggf{\varphi}{R} -\ggf{\varphi^*}{R}\|_{L^2(\Aset, \pi)} 
&\leq \|\Delta (\varphi - \varphi^*)\|_{L^2(\Aset, \pi)} +
\|\langle \nabla (\varphi - \varphi^*), \nabla\log\pi\rangle \|_{L^2(\Aset, \pi)} \\
&\leq \smm{i=1}{d}\|\frac{\partial^{2}(\varphi - \varphi^*)}{\partial x_i^2}\|_{L^2(\Aset, \pi)} + \smm{i=1}{d}\|\frac{\partial(\varphi - \varphi^*)}{\partial x_i}\|_{L^4(\Aset, \pi)}\|\frac{\partial \log \pi}{\partial x_i}\|_{L^4(\Aset, \pi)}\eqsp.
\end{align*}
Since $\pi$ is a probability measure, $\|\cdot\|_{L^2(\Aset, \pi)} \leq \|\cdot\|_{L^\infty(\Aset)}$, and it is easy to see that  
\[
\|\frac{\partial \log \pi}{\partial x_i}\|_{L^4(\Aset, \pi)} \leq \max_\Aset (U'_i(x)) \leq \pg(R)\eqsp,
\]
where the last inequality is due to \Cref{assu:AUF}. Thus, we get
\begin{align*}
\|\ggf{\varphi}{R} -\ggf{\varphi^*}{R}\|_{L^2(\Aset, \pi)} 
&\leq \smm{i=1}{d}\|\frac{\partial(\varphi - \varphi^*)}{\partial x_i^2}\|_{L^\infty(\Aset)} + \pg(R) \smm{i=1}{d} \|\frac{\partial(\varphi - \varphi^*)}{\partial x_i}\|_{L^\infty(\Aset)} \\
&\leq d\max\limits_{|\bgamma| \leq 2} \|D^{\bgamma} (\varphi - \varphi^*)\|_{L^\infty(\Aset)} + d\cdot \pg(R)\max\limits_{|\bgamma| \leq 2} \|D^{\bgamma} (\varphi - \varphi^*)\|_{L^\infty(\Aset)} \\
&= (1+\pg(R))d \|\varphi - \varphi^*\|_{C^2(\Aset)}
\end{align*}
\end{proof}
    
\begin{Lem}
\label{lem:gamma_func_beh}
Assume \Cref{assu:AUF}. Let $f:\rset^d \to \rset$ be a polynomial function with $\pdeg(f)\leq k$. Then 
for $R\geq 1$ it holds that

\[
\|f\|^2_{L^2(\setcomp{\Aset}, \pi)} 
\lesssim
R^{d+2k}e^{-R^2}
\]
up to multiplicative constant $C=C(\mu, d, k)$.
\end{Lem}
\begin{proof}
Recall, that $\mu$-strong convexity of function $U(x)$ is equivalent to 
\[
U(y)
\geq 
U(x)+
\langle
\nabla U(x),y-x
\rangle
+
\frac{\mu}{2}\|x-y\|^2,\quad
\forall x,y\in \rset^d.
\]
It is clear that $x^*$ can be found for strongly convex function such that $\nabla U(x^*) = 0$. This way for $R\geq 1$ we have
\begin{align*}
\|f\|^2_{L^2(\setcomp{\Aset}, \pi)} 
&\lesssim
\intg{\setcomp{\Aset}}{}
f^2(x) e^{-U(x)}dx 
\leq
\intg{\setcomp{\Aset}}{}
f^2(x) e^{-U(x^*)-\mu \|x-x^*\|^2/2}dx 
\\
&\lesssim
\intg{\setcomp{\Aset}}{}
f^2(x) e^{-\|x\|^2}dx 
\lesssim
\intg{R}{\infty}
x^{d+2k} e^{-x^2}dx,
\end{align*}
where last inequality is by polar change of variables. The rest of the proof follows from integration by parts:

\begin{align*}
\intg{R}{\infty}
x^{n} e^{-x^2}dx
&=
\frac{R^{n-1}}{2}e^{-R^2}
+
\frac{n-1}{2}\intg{R}{\infty}
x^{n-2} e^{-x^2}dx= \{ \text{by induction} \}
\\
&\leq 
e^{-R^2}(n-1)!!
\left(\smm{i=0}{\lceil n/2\rceil-1}\frac{R^{n-2i-1}}{(n-2i-1)!!\cdot 2^{i+1}}
+
\intg{R}{\infty}
e^{-x^2}dx\right)
\lesssim
R^{n}e^{-R^2}\,.
\end{align*}
\end{proof}

\section{Appendix}
\label{sec:appendix}
\appendix
\section{Regularity of elliptic PDE solution}
\label{sec:pde_solution_regularity}
We first state the general result on regularity of the solution to a certain partial differential equation. Define the elliptic operator 
\begin{equation}
\label{eq:L_oper}
\textstyle{
\mathcal{L} u =
\sum_{i,j=1}^d a_{ij}(x)\frac{\partial^2 u}{\partial x_i \partial x_j}+\sum_{i=1}^d b_{i}(x)\frac{\partial u}{\partial x_i}+c(x)u}\eqsp,
\end{equation}
and consider a partial differential equation on some domain $B \subseteq \rset^d$:
\begin{equation}\label{eq:dirichletarray}
\mathcal{L} u = f 
\end{equation}
Consider the following assumptions of the coefficients of \eqref{eq:L_oper}, its right-hand side $f(x)$, and the domain $B$:

\begin{myassump}{Reg}
\label{assum:reg}
Let the domain $B = B_0(R) = \{x\in\rset^d:\|x\|_2 < R\}$, $a_{ij}(x)=\indibr{i = j}$, $c(x) = 0$, $b_i = \frac{\partial \log{\pi}}{\partial x_i}, b_i, f \in \H^{s+\varpi}(B_0(R), C_{\ref{assum:reg}} R^{\mu_1})$. 
\end{myassump}
Note that under \Cref{assum:reg} the equation \eqref{eq:dirichletarray} corresponds to the Stein equation \eqref{eq:diff_eq_alt}. That is why we need to study the regularity properties of \eqref{eq:dirichletarray}, namely, we need to bound the norm $\|u\|_{\mathcal{H}^{s+2+\varpi}(B_0(R))}$. The goal of this section is to upper bound the norm 
\[
\|u\|_{\mathcal{H}^{s+2+\varpi}(B_0(R))}\leq K_1 R^{K_2}
\]
of solution to \eqref{eq:dirichletarray} for $R \geq R_0$ under \Cref{assu:AUF}. The constant $K_2$ in the bound above should be \textit{independent} of the problem dimension $d$.
\par 
We begin with introducing multi-indices $\bgamma = (\gamma_1,\dots,\gamma_d) \in \nzeroset^d$ and $|\bgamma| = \sum_{i=1}^{d}\gamma_i$, we define the following norms from \cite[Section~8, Chapter~3]{friedmande}:
\begin{equation}
\label{eq:definitios_pde}
\begin{split}
&d_x = \inf_{z \in \partial B} \|x - z\|\,, \quad  d_{xy} = \min(d_x,d_y), \quad 
H_\varpi (d^k u)(B) = \sup_{x\neq y\in B} d_{xy}^{k+\varpi}\frac{\|u(x)-u(y)\|_\infty}{\|x-y\|_\infty^\varpi}, \\
&|d^k v|_0(B) = \sup_{x\in B} |d_x^k v(x)|\eqsp, \quad 
|u|_m(B) = \sum_{ |\bgamma| =0}^m |d^{|\bgamma|} D^{\bgamma} u|_0(B), \\
&|u|_{m+\varpi}(B) = |u|_m(B) + \sum_{ |\bgamma|=0}^m H_\varpi(d^{|\bgamma|} D^{\bgamma} u)(B)\eqsp.
\end{split}
\end{equation}
We omit the explicit dependence of the quantities above upon the domain $B$ when there is no risk of confusion. We begin with the following non-quantitative result:

\begin{Lem}[\protect{\cite[Theorem~5.6.3]{morrey2008}}]
\label{prop:weak_to_smooth} Let $\mathcal{L}$ be an elliptic operator \eqref{eq:L_oper}, Assume that $a_{ij}(x)=\indibr{i = j}$, $c(x) = 0$, $b_i, f\in \H^{s+\varpi}(B_0(R), C R^{\mu_1})$ for $C,R,\mu_1 > 0, s\in \nzeroset, 0< \varpi \leq 1$. Moreover, assume there exists $u\in W^{2,p}_{loc}(\rset^d)$ satisfying $\mathcal{L}u=f$, $d <p$. Then $u\in C^{s+2+\varpi}(B_0(R))$.
\end{Lem}
\begin{proof}
We have that $u\in W^{2,p}(B)$ for every $B \subset \rset^d$. By assumptions on coefficients and $f$, \cite[Theorem~5.6.3, (b)]{morrey2008} implies that that $u\in C^{s+2+\varpi}(B_0(R))$.
\end{proof}

Note that \Cref{prop:weak_to_smooth} does not automatically imply even that $\|u\|_{\mathcal{H}^{s+2+\varpi}(B_0(R))} < \infty$. To bound $\|u\|_{\mathcal{H}^{s+2+\varpi}(B_0(R))}$, we proceed with the interior Schauder estimates introduced in \cite[Theorem~5, Chapter~3]{friedmande}. Under the assumptions of \Cref{prop:weak_to_smooth}, it implies that 
\begin{equation}
\label{eq:base_interior_bound}
|u|_{2+\varpi}\leq K_{R} (|u|_0 + |d^2 f|_\varpi)\,,
\end{equation}
where the constant $K_{R}$ possibly depends upon $R$. Moreover, due to the definitions \eqref{eq:definitios_pde}, one can check that 
\[
\|u\|_{\mathcal{H}^{s+\varpi}(B_0(R))} \lesssim |u|_{s+\varpi}(B_0(R+1)),\quad 
|u|_{s+\varpi}(B_0(R)) \lesssim R^{s+\varpi} \|u\|_{\mathcal{H}^{s+\varpi}(B_0(R))}\eqsp,
\]
where inequalities hold up to constant depending $s,\varpi, d$.
Thus, it is enough to bound $|u|_{s+\varpi}$ on $B_0(R+1)$. Tracing the proof of \cite[Theorem~5, Chapter~3]{friedmande} we state the following result: 
\begin{Lem}
\label{prop:const_bound} Assume \Cref{assum:reg}. Then the solution $u(x)$ of \Cref{eq:dirichletarray} satisfies
\[
|u|_{2+\varpi} \leq h(\max(1,C_{\ref{assum:reg}} R^{\mu_1}), R ,|u|_0, |d^2 f|_\varpi)\eqsp,
\]
and $h$ is a polynomial function with $\pdeg(h)$ being independent of $d$.
\end{Lem}
\begin{Rem}
    The bound \eqref{eq:base_interior_bound} is proved in \cite[Th.~1, Ch.~4]{friedmande}. In formulation of theorem it is given, that $K$ is independent of the diameter $R$. Using \eqref{eq:base_interior_bound}, the only thing to trace is the polynomial scaling of $K_R$.
\end{Rem}

Our next step is to get the bounds for $|u|_{s+\varpi}$ instead of just $|u|_{2+\varpi}$. Usually, we have so-called apriori bounds on $|u|_0$ for solutions of PDEs, which are in general exponential w.r.t. $R$. One can find an overview in \cite[Chapter~9.1]{gilbarg1977elliptic}. We need more accurate approach which gives us the bound only for some particular solution. Namely, we use the following result from \cite{bogachev2018poisson}:

\begin{Lem}[\protect{\cite[Example~4.9, Part~4]{bogachev2018poisson}}]
\label{lem:good_bound} Consider the equation \eqref{eq:dirichletarray} with operator $\mathcal{L}$ defined in \eqref{eq:L_oper}. Let $a_{ij} \in W_{loc}^{1,p}(\rset^d)$, $b_i\in L_{loc}^p(\rset^d)$, $p > d$, $(a_{ij}(x))_{i,j=1}^{d}$ is symmetric, positive definite, uniformly bounded, $b$ satisfies
\[
\langle b, x\rangle \leq C_{\ref{lem:good_bound}, 1}-C_{\ref{lem:good_bound}, 2}\|x\|^\gamma\eqsp,
\]
function $f$ satisfies 
\begin{equation*}
    \|f(x)\|\leq C_{\ref{lem:good_bound}, 2}(1+\|x\|)^{2\alpha+\gamma-1}\eqsp,
\end{equation*}
and 
\[
\int_{\rset^d} f d\mu = 0\eqsp, 
\]
where $\mu$ satisfies
\[
\int_{\rset^d}\mathcal{L} u(x) \, d\mu(x) = 0,\quad \forall u\in C_0^\infty(\rset^d)\eqsp.
\]
Then there exists a unique function $u\in W^{2,p}_{loc}(\rset^d)$ satisfying $\mathcal{L}u= f$ such that
\[
\|u(x)\|\leq C_{\ref{lem:good_bound}, 3}(1+\|x\|)^{2\alpha}
\]
where $C_{\ref{lem:good_bound}, 1}, C_{\ref{lem:good_bound}, 2}, C_{\ref{lem:good_bound}, 3}$ are some constants, independent of $x$.
\end{Lem}

Now we generalize the result of 
\Cref{prop:const_bound} for higher-order derivatives.  

\begin{Lem}
\label{lem:sol_holder_const} 
Assume \Cref{assum:reg} and \Cref{assu:AUF}. Then the solution $u(x)$ of \eqref{eq:dirichletarray} belongs to $\mathcal{H}^{\beta+3}(D, C_{\ref{lem:sol_holder_const}, 2})$ with the constant
\[
C_{\ref{lem:sol_holder_const}, 1} \leq h(\max(1,C_{\ref{assum:reg}} R^{\mu_1}), R)\eqsp.
\]
Here $h$ is a polynomial function with $\pdeg(h) = \hdeg$ independent of $d$.
\end{Lem}

\begin{proof}
We start with \Cref{lem:good_bound}. It gives us function $u\in W^{2,p}_{loc}(\rset^d)$ satisfying equation \eqref{eq:dirichletarray} and growing not faster than polynomially. Next we check the conditions of \Cref{lem:good_bound}. Namely, under \Cref{assu:AUF} we get
\[
\langle b(x), x\rangle = \langle \nabla\log\pi(x), x\rangle = -\langle \nabla U(x), x\rangle
\]
Without loss of generality, let $\nabla U(0)=0$. Thus by definition of strong convexity
\[
-\langle \nabla U(x), x\rangle \leq -\|x\|_2^2\eqsp.
\]
Noe that we can set $p>d$. Measure $\pi$ satisfies the conditions of \Cref{lem:good_bound} from integratio by parts, as it is discussed in \Cref{sec: setup}.
    Under \Cref{assu:AUF}, \Cref{prop:weak_to_smooth} ensures that $u \in C^{2+s+\varpi}$, and \Cref{prop:const_bound} yields 
    \[
    \|u\|_{\H^{2+\varpi}(B_0(R))}\leq |u|_{2+\varpi}(B_0(R+1))\eqsp.
    \]
    \par 
    Next, it is left to generalise \Cref{prop:const_bound} to $\|u\|_{\H^{s+2+\varpi}(B_0(R))}$. Under smoothness assumptions in \Cref{assu:AUF} we get the polynomial dependence on $C_{\ref{lem:sol_holder_const}, 1}$ and $R$ and degree of this polynomial dependence is not dependent on $d$ in the same way as in \cite[Th.~10, Ch.~3]{friedmande}.
\end{proof}
    
\section{Rate of convergence for the variance excess}
\label{sec:excess_belomestny_statcomp}
    
    \begin{myassump}{B}
    \label{assu:B}
    There exist $B>0$ such that \(\sup_{h\in\H}{|h|_{\infty}}\leq B\) with \({|h|_{\infty}}=\sup_{x\in \Xset}|h(x)|\)
    \end{myassump}

    
    \begin{Th}[\protect{\cite[Theorem~3]{belomestny_variance_2020_esvm}}]\label{th:main_slow}
    Assume \Cref{assu:ge} and \Cref{assu:B} + \Cref{assu:br}. Set \(b_n=2(\log(\fraca{1}{\rho}))^{-1}\log(n)\) and take
    	$\eps = \gamma_{\ltwo(\pi)}(\mathcal{H},n)$. Let
     \[
     \widehat{h}_{\eps} \in \argmin_{h \in \mathcal{H}_{\eps}} \SV{h}\,.
     \]
     Then for any \(\delta \in (0,1),\) there exists \(n_0=n_0(\delta)>0\) such that for any $n \geq n_0$ and $x_0 \in \Cset$ with $\P_{x_0}-$probability at least $1-\delta,$ it holds
    	\begin{equation*}\label{eq:th1}
    		\AV{\hh{\eps}} - \inf_{h\in\H}\AV{h}
    		\lesssim
    		 C_{\ref{th:main_slow}, 1} \log(n) \gamma_{\ltwo(\pi)}(\mathcal{H},n)
    		+ C_{\ref{th:main_slow}, 2}\frac{\log(n)\log(\fraca{1}{\delta})}{\sqrt{n}},
    	\end{equation*}	
    	where \(\lesssim\) stands for inequality up to an absolute constant,
    	\begin{equation*}
    		C_{\ref{th:main_slow}, 1} = \frac{\beta B^2}{\log(\fraca{1}{\rho})}\,,\,
    		C_{\ref{th:main_slow}, 2} =
    		\frac{\varsigma^{1/2} (\pi(\driftfunc) + \driftfunc(x_0))}{(1-\rho)^{1/2}\log(\fraca{1}{\rho})}
    		\biggl(\beta B^2 + \sup_{h \in \H}\VnormFunc[2]{h}{\driftfunc^{1/2}}\biggr)\,,\, \beta=\frac{\varsigma \ls}{1- \rho}\biggl(\frac1{\log u}+\frac{J \varsigma \ls }{1-\rho}\biggr)\eqsp.
    	\end{equation*}
    \end{Th}
    
\section{Approximation with neural networks}
\label{sec:neural_networks_recu}
In \cite[Theorem~1]{belomestny2022simultaneous} the authors study the approximation properties of neural networks with ReQU (rectified quadratic unit) activation function $\sigma^{\mathsf{ReQU}}(x) = (x \vee 0)^2$. Recall that the application of operator \eqref{eq:steincv} requires at least $2$ times continuously differentiable function $\varphi$. Thus, we need the counterpart of \cite[Theorem~1]{belomestny2022simultaneous}, but with rectified cubic units (ReCU) as activation functions. Below we provide the corresponding modification. Albeit the proof technique is similar to \cite[Theorem~1]{belomestny2022simultaneous}, our result has some new desirable features. First, the number of layers of the approximation network is independent of the H\"older constant of target function. Moreover, we need a precise non-asymptotic bound on the number of layers in the constructed network. The respective result is stated below in \Cref{prop:recu_approx_const}. 
\par 
Our proof construction is based on approximating the target function in an appropriate spline basis (using B-splines). For completeness we give the definition of B-spline. For given $q, K\in\nset$ consider a vector $a=(a_1,\dots, a_{2q +K+1})$ such that
\begin{align*}
    &a_1 =\dots= a_{q+1}=0,\\
    &a_{q+1+j} = j/K,\ 1\leq j\leq K-1,\\
    &a_{q+K+1} = \dots = a_{2q+K+1} = 1\,.
\end{align*}
For any $j\in\{1,\dots, 2q+K\}$ let
\[
B_j^{0,K}=
\begin{cases}
\frac{1}{a_{j+1}-a_j },\ \text{if } a_j< a_{j+1},
a_j\leq x < a_{j+1},\\
0,\ \text{otherwise}\,.
\end{cases}
\]
Then for $m\in\{1,\dots, q\}$, $j\in\{1,\dots, 2q+K-m\}$ define
\begin{equation}
\label{eq:b_spline_def}
B_j^{m, K}(x)=
\begin{cases}
\frac{(x-a_j)B_j^{m-1, K}(x) + (a_{j+m+1}-x)B_{j+1}^{m-1, K}(x) }{a_{j+m+1}-a_j },\ \text{if } a_j< a_{j+m+1},
a_j\leq x < a_{j+m+1},\\
0,\ \text{otherwise}\,.
\end{cases}
\end{equation}
Now we state the main results of this section.

\begin{Th}\label{prop:recu_approx_const}
Fix  $\beta > 3$ and $p, d \in \mathbb N$.
	Then, for any $f : [0, 1]^d \rightarrow \R^p$, $f \in \H^\beta([0,1]^d, H)$ and any integer $K \geq 3$, there exists a neural network $h_f : [0, 1]^d \rightarrow \R^p$ with ReCU activation functions  such that it has at most $8\cdot(4+2(\lbeta-2)+\lceil\log_2 d\rceil + 3)$ layers,  at most $O(d\vee p \vee (K+\lbeta)^{d}(2\lbeta+1)^{d}9^{d(\beta-1)}H)$ neurons in each layer and $O(d\lbeta (K+\lbeta)^{d}(2\lbeta+1)^{d}9^{d(\beta-1)}H)$ non-zero weights taking their values in $[-1, 1]$, such that
    \[
    \left\|f - h_f\right\|_{\H^\ell([0, 1]^d)} \leq \frac{C_{\ref{prop:recu_approx_const}, 1}^{\beta d} H \beta^\ell}{K^{\beta-\ell}}
	    \quad
	    \text{for all $\ell \in \{0, \ldots, \lfloor\beta\rfloor\}$,}
    \]
where $C_{\ref{prop:recu_approx_const}, 1}$ is an absolute constant.
\end{Th}

\begin{proof}
The proof follows the main steps of \cite[Theorem~2]{belomestny2022simultaneous}. The main differences with the mentioned paper are as follows. First, we need to estimate number of additional layers required due to using ReCU activations instead of ReQU. Second, we modify the construction of multiplication by a number given in \cite[Lemma~4]{belomestny2022simultaneous} in order to implement it with a fixed number of layers, with possible ncrease of their width.
\par 
Recall that $\sigma(x) = (0\vee x)^3$. Note that
\[
\sigma(x) + \sigma(-x) = x^3
\implies
x^3 =\begin{pmatrix}
1&1
\end{pmatrix}\circ
\sigma\circ
\begin{pmatrix}
1\\
-1
\end{pmatrix}\circ x \eqsp,
\]
thus it takes 2 layers to get $x^3$. From now on we use numbers above arrows to denote number of layers required to implement the corresponding transformation. Thus we can construct mappings
\[
x\xrightarrow{2} ((x+1)^3, (x-1)^3, x^3)\xrightarrow{1} 
\frac{(x+1)^3 + (x-1)^3-2x^3}{6}
=
x,
\]
\[
x\xrightarrow{2} ((x+1)^3, (x-1)^3)\xrightarrow{1} 
\frac{(x+1)^3 - (x-1)^3-2}{6}
=
x^2,
\]
which both requires a network with $3$ layers. Thus it takes 4 layers to implement the product
\begin{multline*}
(x_1, x_2)\xrightarrow{2}
((x_1+x_2+1)^3, (x_1+x_2-1)^3, (x_1-x_2+1)^3, (x_1-x_2-1)^3)
\\
\xrightarrow{1}
((x_1+x_2)^2, (x_1-x_2)^2)
\xrightarrow{1}\frac{(x_1+x_2)^2 - (x_1-x_2)^2}{4}= x_1x_2\,,
\end{multline*}
and $8$ layers to implement 
\begin{equation*}
x\xrightarrow{1} \underbrace{(1, 1, \dots, 1)}_{\text{K times}}
\xrightarrow{2}
((K+1)^3, (K-1)^3)\xrightarrow{1}(K^2, K^2)
\xrightarrow{4} K^4\,.
\end{equation*}
Our next step is to construct a mapping 
\[
x\mapsto 
(x, K, B_1^{3, K}(x),
B_2^{3, K}(x),
\dots, 
B_{K+2q-3}^{3, K}(x)
),
\]
where $B_{j}^{m,K}$ is defined in \eqref{eq:b_spline_def}, and $q, K \geq 3$, $x\in [0;1]$. Proceeding as in \cite[Lemma~2]{belomestny2022simultaneous}, we express $B_j^{m, K}$ using $x_+^3$. Note that with ReCU activation we can only express the B-splines with indices $m, K \geq 3$.
\par 
Finally, we implement mapping $x\to Mx$ in a more compact way with respect to number of layers. We concatenate maps
\[
x\xrightarrow{4} x,\;
x\xrightarrow{2} \underbrace{(1, \dots, 1)}_{M\text{ times}}
\xrightarrow{2} M
\]
and obtain
\[
(M,x)\xrightarrow{4} Mx,\qquad 
x \xrightarrow{8} Mx\,.
\]
With this multiplication construction, we trace the proof of \cite[Theorem~2]{belomestny2022simultaneous}, and get a neural network of depth at most 
\[
L \leq 8\cdot(4+2(\lbeta-2)+\lceil\log_2 d\rceil + 2)+ 8
\]
and width $O(d\vee p \vee (K+\lbeta)^{d}(2\lbeta+1)^{d}9^{d(\beta-1)}H)$. The rest of the proof remains unchanged.
\end{proof}

Recall definition of neural network \eqref{eq:nn}. Following \cite[Lemma 5]{schmidt-hieber2020}, we additionally introduce the following notations. For $k \in \{1,\dots,L+1\}$, $i \leq k$, $k \leq j \leq L$, we define functions $\BBB_{k,i}(x): \R^{p_{i-1}} \rightarrow \R^{p_{k}}$ and $\AAA_{j,k}(x): \R^{p_{k-1}} \rightarrow \R^{p_{j+1}}$ as follows
\begin{equation}
\label{eq:A_plus_A_minus_def}
\begin{split}
\BBB_{k,i}(x) &= \requ_{v_{k}} \circ W_{k-1} \circ \requ_{v_{k-1}} \dots \circ \requ_{v_i} \circ W_{i-1} \circ x; \\
\AAA_{j,k}(x) &= W_{j}  \circ \requ_{v_{j}} \dots \circ W_{k} \circ \requ_{v_{k}} \circ W_{k-1} \circ x\,.
\end{split}
\end{equation}
Set by convention $\AAA_{L,L+2}(x) = \BBB_{0,1}(x) = x$. For notation simplicity, we write $\AAA_{j}$ instead of $\AAA_{L,j}$, and $\BBB_{j}$ instead of $\BBB_{j,1}$. Note that with this notation $f(x) = \AAA_{1}(x)$.
\par
Let us introduce the functions 
\begin{equation}
\label{eq:f_one_f_two_def}
\begin{split}
\fone(x) &= \Wone_L \circ \requ_{\vone_{L}} \circ \Wone_{L-1} \circ \requ_{\vone_{L-1}} \circ \dots \circ \requ_{\vone_1} \circ \Wone_0 \circ x\,, \\
\ftwo(x) &= \Wtwo_L \circ \requ_{\vtwo_{L}} \circ \Wtwo_{L-1} \circ \requ_{\vtwo_{L-1}} \circ \dots \circ \requ_{\vtwo_1} \circ \Wtwo_0 \circ x\,,
\end{split}
\end{equation}
where the parameters $\Wone_{i}, \Wtwo_{i}, \vone_{i}, \vtwo_{i}$ satisfy
\[
\norm{\Wone_{i} - \Wtwo_{i}}_{\infty} \leq \eps, \quad \norm{\vone_{i} - \vtwo_{i}}_{\infty} \leq \eps\,, \quad \text{for all } i \in \{0,\dots,L\}\,.
\]

\begin{Prop}
\label{Lem:h2_norm_bound}
Consider neural networks of class $\NN(L, \A, \snzw)$ with architecture $\A = (p_0,\dots,p_{L},p_{L+1})$ and at most $\snzw$ non-zero weights satisfying $\norm{W_{i}}_{\infty} \leq 1, \norm{v_i}_{\infty} \leq 1$, $i \in \{0,\dots,L\}$. Then for any $\eps > 0$
\begin{equation*}
\label{eq:cov_number_bound}
\mathcal{N}(\NN(L, \A, \snzw), \|\cdot\|_{C^2([0,1]^{p_0})}, \eps) 
\leq \left(2\eps^{-1} L^2(L+1)(L+2)3^{3L+4}
V^{L3^{L+5}+1}\right)^{\snzw+1}\eqsp,
\end{equation*}
where we have denoted
\begin{equation}
\label{eq:V_def}
V = \prod_{\ell=0}^{L}(p_{\ell}+1)\,.
\end{equation}
\end{Prop}
\begin{proof}
We follow \cite[Lemma 5]{schmidt-hieber2020}. Fix arbitrary $\eps > 0$ and let $\fone(x), \ftwo(x) \in \NN(L, \A, \snzw)$ be neural networks with weights $\norm{\Wone_{i} - \Wtwo_{i}}_{\infty} \leq \eps, \norm{\vone_{i} - \vtwo_{i}}_{\infty} \leq \eps$, $i \in \{0,\dots,L\}$. Then \Cref{Lem:inf_norm_bound} and \Cref{Lem:inf_norm_bound_derivative} imply 
\begin{align*}
\norm{\fone(x) -  \ftwo(x)}_{C^2([0;1]^{p_0})} \leq \eps L^2(L+1)(L+2)3^{3L+4}
V^{L3^{L+5}}\,.
\end{align*}
Note that for functions in class $\NN(L, \A, s)$ number of parameters is bounded by $T=\sum_{\ell=0}^{L}(p_{\ell}+1)p_{\ell+1} \leq (L+1)2^{-L}\prod_{\ell=0}^{L+1}(p_\ell+1)\leq V$, hence, there are at most $\binom{T}{\snzw}\leq V^{\snzw}$ possible combinations to allocate $\snzw$ non-zero weights. Since all parameters are bounded in absolute value by one, we can discretize the non-zero parameters with grid step $\delta/(2L^2(L+1)(L+2)3^{3L+4}
V^{L3^{L+5}})$ and obtain
\begin{align*}
\mathcal{N}(\NN(L, \A, \snzw), \|\cdot\|_{C^2([0,1]^{p_0})}, \delta) 
&\leq \sum_{\ell = 1}^{\snzw}\left(2\delta^{-1} L^2(L+1)(L+2)3^{3L+4}
V^{L3^{L+5}}\right)^{\ell}\cdot V^\ell  \\ 
&\leq \left(2\delta^{-1} L^2(L+1)(L+2)3^{3L+4}
V^{L3^{L+5}+1}\right)^{\snzw+1}\,.
\end{align*}
The bound for $C^1$ and $C^0=L_\infty$ metrics can be obtained in a similar manner.
\end{proof}

\begin{Prop}
\label{lem:infty_norm_hessian}
Let $x \in \R^{d}$, $\norm{x}_{\infty} \leq \Kxnorm$, $\nabla^2 f(x)$ be from \eqref{eq:sec_grad}. Then for $k, i \in \{1,\dots,L\}, \, k \geq i$ 
\begin{align*}
\norm{\nabla^2 f(x)}_{\infty}  \leq
2^5 \cdot 144^{L}\cdot L
\cdot 
    \biggl\{\prod_{\ell=0}^{L}(p_{\ell}+1)^{4(L3^{L}+1)+1}\biggr\} (\Kxnorm \vee 1)^{4\cdot 3^{L}}\,. 
\end{align*}
\end{Prop}
\begin{proof}
Recall the definition:
\begin{align*}
(\nabla^2 f(x))^T \notag &= \biggl(
\smm{i=1}{L}\smm{j=1}{p_{L-i+1}}
2\cdot 3^{L}\cdot 
\biggl[W_{L} \bigl(
\prod_{\ell=1}^{i-1} 
\left\{\operatorname{diag}\left[(\AAA_{L-\ell,1}(x) + v_{L-\ell+1} \vee 0)^2\right] W_{L-\ell}\right\}
\bigr) \notag\\
&\times
e_j^T \operatorname{diag}\left[\AAA_{L-i,1}(x) + v_{L-i+1} \vee 0\right]\cdot 
J(\AAA_{L-i,1}(x))\biggr]^T \cdot e_j e_j^T
 \notag\\
&\times
W_{L-i}
\bigl(
\prod_{\ell=i+1}^{L} \left\{\operatorname{diag}\left[(\AAA_{L-\ell,1}(x) + v_{L-\ell+1} \vee 0)^2\right] W_{L-\ell}\right\}
\bigr)
\biggr)\,. 
\end{align*}
Hence,
\begin{align*}
&\norm{\nabla^2 f(x)}_{\infty} \notag \leq \biggl(
\smm{i=1}{L}\smm{j=1}{p_{L-i+1}}
2\cdot 3^{L}\cdot 
\biggl[(p_{L}+1)\bigl(
\prod_{\ell=1}^{i-1}
(p_{L-\ell}+1)^2(\norm{\AAA_{L-\ell,1}}_{\infty}+1)^2
\bigr) \notag\\
&\cdot
(p_{L-i}+1)^{2}(\norm{\AAA_{L-i,1}}_{\infty}+1)
\norm{J(\AAA_{L-i,1})}_\infty \bigl(
\prod_{\ell=i+1}^{L}
(p_{L-\ell}+1)^2(\norm{\AAA_{L-\ell,1}}_{\infty}+1)^2
\bigr)
\biggr)\\
&\leq 
2\cdot 3^L\cdot L\cdot(\max_{0\leq i\leq L}p_i)\cdot (p_{L}+1)
\bigl(\prod_{\ell=1}^{L}
(p_{L-\ell}+1)^2(\norm{\AAA_{L-\ell,1}}_{\infty}+1)^2
\bigr)(\max_{0\leq i\leq L}\norm{J(\AAA_{L-i,1})}_\infty)
\,. 
\end{align*}
To bound $\norm{J(\AAA_{L-i,1})}_\infty$ recall the definition
\begin{equation*}
J(\AAA_{L})^T=
(\nabla f(x))^T = 3^{L} W_{L} \prod_{\ell=1}^{L} \left\{\operatorname{diag}\left[(\AAA_{L-\ell,1}(x) + v_{L-\ell+1} \vee 0)^2\right] W_{L-\ell}\right\}\,.
\end{equation*}
Thus,
\begin{equation*}
\norm{J(\AAA_{L})}_\infty \leq
3^{L} (p_{L}+1)
\bigl(\prod_{\ell=1}^{L}
(p_{L-\ell}+1)^2(\norm{\AAA_{L-\ell,1}}_{\infty}+1)^2
\bigr)\,.
\end{equation*}
Using the bound of \Cref{lem:infty_norm_a_plus}, we get 
\begin{equation*}
    \norm{\AAA_{L-\ell,1}}_{\infty}+1\leq 
    2(p_L+1)\biggl\{\prod_{i=1}^{L-\ell}(p_{L-\ell-i}+1)^{3^{i}}\biggr\} (\Kxnorm \vee 1)^{3^{L-\ell}}\eqsp,
\end{equation*}
and
\begin{align*}
    \prod_{\ell=1}^{L}
    (\norm{\AAA_{L-\ell,1}}_{\infty}+1)
    &= 2^L(p_L+1)^L
    \prod_{\ell=1}^{L}
    \biggl(\biggl\{\prod_{i=1}^{L-\ell}(p_{L-\ell-i}+1)^{3^{i}}\biggr\} (\Kxnorm \vee 1)^{3^{L-\ell}}
    \biggr)\\
    &\leq 2^L(p_L+1)^L
    \biggl\{\prod_{\ell=0}^{L}(p_{\ell}+1)^{L3^{L}}\biggr\} (\Kxnorm \vee 1)^{3^{L}}\eqsp.
\end{align*}
Combining the bounds gives the result.
\end{proof}

\begin{Lem}
\label{lem:infty_norm_a_plus}
Let $x \in \R^{d}$, $\norm{x}_{\infty} \leq \Kxnorm$. Then for $k, i \in \{1,\dots,L\}, \, k \geq i$ 
\begin{equation}
\label{eq:A_plus_infty_bound}
\norm{\BBB_{k,i}(x)}_{\infty} \leq \biggl\{\prod_{\ell=1}^{k-i+1}(p_{k-\ell}+1)^{3^{\ell}}\biggr\} (\Kxnorm \vee 1)^{3^{k-i+1}}\,,
\end{equation}
Moreover, function $\AAA_{j,k}(x)$ is Lipshitz for $x,y \in \R^d: \norm{x}_{\infty} \leq \Kxnorm, \norm{y}_{\infty} \leq \Kxnorm$, that is,
\begin{equation}
\label{eq:A_minus_lip_bound}
\norm{\AAA_{j,k}(x) - \AAA_{j,k}(y)}_{\infty} \leq 3^{j-k+1} \prod_{\ell=0}^{j-k+1}(p_{j-\ell}+1)^{2\cdot 3^{\ell}}(\Kxnorm \vee 1)^{2\cdot 3^{j-k+1}}\norm{x-y}_{\infty}\,. 
\end{equation}
\end{Lem}
\begin{proof}
\eqref{eq:A_plus_infty_bound} follows from an easy induction in $k$. Indeed, if $k = i$, $\BBB_{i,i} = \requ_{v^{(i)}} \circ W_{i-1}x$, and 
\[
\norm{\BBB_{i,i}}_{\infty} \leq (\Kxnorm p_{i-1}+1)^3 \leq (p_{i-1}+1)^3(\Kxnorm \vee 1)^{3}\,.
\]
Using $\norm{\BBB_{k,i}}_{\infty} \leq (\norm{\BBB_{k-1,i}}_{\infty} p_{k-1}+1)^3$ completes the proof. 
\par 
To prove \eqref{eq:A_minus_lip_bound}, we use an induction in $j$. Assume that \eqref{eq:A_minus_lip_bound} holds for any $k \in \{1,\dots,L+1\}$ and $j - 1 \geq k$. Then
\begin{align*}
&\norm{\AAA_{j,k}(x) - \AAA_{j,k}(y)}_{\infty} 
\leq p_j \norm{\BBB_{j,k}(x) - \BBB_{j,k}(y)}_{\infty} \\
&\leq 3p_{j}\norm{\AAA_{j-1,k}(x) - \AAA_{j-1,k}(y)}_{\infty} \bigl(\norm{\AAA_{j-1,k}(x)}_{\infty} \vee \norm{\AAA_{j-1,k}(y)}_{\infty}\bigr)^2 \\
&\leq 3(p_{j}+1)(p_{j-1}+1) \prod_{\ell=1}^{j-k}(p_{j-\ell-1}+1)^{2\cdot 3^{\ell}}(\Kxnorm \vee 1)^{2\cdot 3^{j-k}} \norm{\AAA_{j-1,k}(x) - \AAA_{j-1,k}(y)}_{\infty}\,,
\end{align*}
and the statement follows from the elementary bound 
\begin{align*}
\norm{\AAA_{k-1,k}(x) - \AAA_{k-1,k}(y)}_{\infty} \leq p_{k-1}\norm{x-y}_{\infty} \leq p_{k-1}\norm{x-y}_{\infty} (\Kxnorm \vee 1)\,.
\end{align*}
\end{proof}

\begin{Lem}
\label{Lem:inf_norm_bound}
For the neural networks $\fone(x)$ and $\ftwo(x)$ defined in \eqref{eq:f_one_f_two_def}, $\|x\|_\infty \leq 1$ it holds
\begin{equation}
\label{eq:inf_norm_bound}
\norm{\fone(x) - \ftwo(x)}_{\infty} \leq \eps (L+1) 3^{L} \prod_{\ell=0}^{L}(p_{\ell}+1)^{3^{L+2}}\,.
\end{equation}
\end{Lem}

\begin{proof}
Denote by $\AAA_{j}^{(1)}, \BBB_{j}^{(1)}, \AAA_{j}^{(2)}, \BBB_{j}^{(2)}$ the corresponding functions in \eqref{eq:A_plus_A_minus_def}. Following \cite[Lemma 5]{schmidt-hieber2020}, we write
\[
\norm{\fone(x) - \ftwo(x)}_{\infty}
\leq \sum_{k=1}^{L+1}\norm{\AAA^{(1)}_{k+1} \circ \requ_{\vone_{k}} \circ \Wone_{k-1} \circ \BBB_{k-1}^{(2)}(x) - \AAA^{(1)}_{k+1} \circ \requ_{\vtwo_{k}} \circ \Wtwo_{k-1} \circ \BBB_{k-1}^{(2)}(x)}_{\infty}\,.
\]
Due to \Cref{lem:infty_norm_a_plus}, functions $\AAA^{(1)}_{k+1} = \AAA^{(1)}_{L,k+1}$ are Lipshitz, and
\begin{multline*}
\norm{\fone(x) - \ftwo(x)}_{\infty} \leq \sum_{k=1}^{L+1}3^{L-k}\prod_{\ell=0}^{L-k}\bigl(p_{L-\ell} + 1\bigr)^{2\cdot 3^{\ell}}\left(\norm{\BBB^{(2)}_{k}(x)}_{\infty} \vee 1\right)^{2\cdot 3^{L-k}}  \times 
\\
\quad \norm{\requ_{\vone_{k}} \circ \Wone_{k-1} \circ \BBB^{(2)}_{k-1}(x) - \requ_{\vtwo_{k}} \circ \Wtwo_{k-1} \circ \BBB^{(2)}_{k-1}(x)}_{\infty}.
\end{multline*}
Note that
\begin{align*}
&\norm{\requ_{\vone_{k}} \circ \Wone_{k-1} \circ \BBB^{(2)}_{k-1}(x) - \requ_{\vtwo_{k}} \circ \Wtwo_{k-1} \circ \BBB^{(2)}_{k-1}(x)}_{\infty}\\ 
&\leq \eps (p_{k-1}+1)\left(\norm{\BBB^{(2)}_{k-1}(x)}_{\infty} \vee 1\right) \times
\biggl(\norm{\Wone_{k-1} \circ \BBB^{(2)}_{k-1}(x)}_{\infty}^2 \\
& + \norm{\Wone_{k-1} \circ \BBB^{(2)}_{k-1}(x)}_{\infty}
\times \norm{\Wtwo_{k-1} \circ \BBB^{(2)}_{k-1}(x)}_{\infty} + \norm{\Wtwo_{k-1} \circ \BBB^{(2)}_{k-1}(x)}_{\infty}^2\biggr) \\
&\leq 3\eps (p_{k-1}+1)^3\left(\norm{\BBB^{(2)}_{k-1}(x)}_{\infty} \vee 1\right)^3\,.
\end{align*}
First inequality is due to difference of $n$-th powers formula and Lipschitz property of ReLU function.
Combining the previous bounds and \eqref{eq:A_plus_infty_bound} yields
\begin{align*}
\norm{\fone(x) - \ftwo(x)}_{\infty} 
&\leq \eps \sum_{k=1}^{L+1}3^{L-k+1}\left\{\prod_{\ell=0}^{L-k}\bigl(p_{L-\ell} + 1\bigr)^{2\cdot 3^{\ell}}\right\} \bigl(p_{k-1} + 1\bigr)^3 \left(\norm{\BBB^{(2)}_{k}(x)}_{\infty}^2 \vee 1\right)^{2\cdot 3^{L-k}} \\
&\leq \eps (L+1) 3^{L} \prod_{\ell=0}^{L}(p_{\ell}+1)^{3^{L+2}}\,.
\end{align*}
and the statement follows.
\end{proof}

\begin{Lem}
\label{Lem:first_sec_derivative}
Let $f(x): \R^{p_0} \rightarrow \R^{p_{L+1}}$ be from \eqref{eq:f_one_f_two_def}. Then for $f(x): \R^{p_0} \rightarrow \R^{p_{L+1}}$ it holds that
\begin{equation}
\label{eq:first_grad}
(\nabla f(x))^T = 3^{L} W_{L} \prod_{\ell=1}^{L} \left\{\operatorname{diag}\left[(\AAA_{L-\ell,1}(x) + v_{L-\ell+1} \vee 0)^2\right] W_{L-\ell}\right\}\,.
\end{equation}
and for $f(x): \R^{p_0} \rightarrow \R$ additionally we get

\begin{align}
\label{eq:sec_grad}
(&\nabla^2 f(x))^T \notag = \biggl(
\smm{i=1}{L}\smm{j=1}{p_{L-i+1}}
2\cdot 3^{L}\cdot 
\biggl[W_{L} \bigl(
\prod_{\ell=1}^{i-1} 
\left\{\operatorname{diag}\left[(\AAA_{L-\ell,1}(x) + v_{L-\ell+1} \vee 0)^2\right] W_{L-\ell}\right\}
\bigr) \notag\\
&\times
e_j^T \operatorname{diag}\left[\AAA_{L-i,1}(x) + v_{L-i+1} \vee 0\right]\cdot 
J(\AAA_{L-i,1}(x))\biggr]^T \cdot e_j e_j^T
 \notag\\
&\times
W_{L-i}
\bigl(
\prod_{\ell=i+1}^{L} \left\{\operatorname{diag}\left[(\AAA_{L-\ell,1}(x) + v_{L-\ell+1} \vee 0)^2\right] W_{L-\ell}\right\}
\bigr)
\biggr)\,. 
\end{align}
\end{Lem}
\begin{proof}
We use matrix differentiation technique. Recall, that for one dimensional outputs it holds
\[
df = \langle \nabla f(x), dx\rangle,\ 
d^2 f = \langle \nabla^2 f(x) dx_2, dx_1\rangle
\]
and for vector input vector output situation we have
$df = J_f dx$, where $J_f$ - Jacobian matrix of function $f$. Thus to find $\nabla^2 f(x)$ we evaluate $d(df) = d(\langle \nabla f(x), h_1\rangle)$ as if difference argument $h_1$ is fixed. First observe, that for $g(x): \R^n\to \R^n$ it holds that
\[
J\bigl(\sigma^{\mathsf{ReQU}}(g(x))\bigr) =
\operatorname{diag}\bigl(2\sigma^{\mathsf{ReLU}}(g(x))\bigr)\cdot J(g(x)),
\]
\[
J\bigl(\sigma^{\mathsf{ReCU}}(g(x))\bigr)
=
\operatorname{diag}\bigl(3\sigma^{\mathsf{ReQU}}(g(x))\bigr)\cdot J( g(x)).
\]
From this it easily follows that

\[
(\nabla f(x))^T = 3^{L} W_{L} \prod_{\ell=1}^{L} \left\{\operatorname{diag}\left[(\AAA_{L-\ell,1}(x) + v_{L-\ell+1} \vee 0)^2\right] W_{L-\ell}\right\}\,.
\]
where $(\AAA_{L-\ell,1}(x) + v_{L-\ell+1} \vee 0)^2$ stands for element-wise raise in power. 
Let $
S_\ell = (\AAA_{L-\ell,1}(x) + v_{L-\ell+1} \vee 0)^2
$. Then
\begin{multline*}
\langle d(\nabla f), h_1\rangle=
d\biggl(3^{L} W_{L} \prod_{\ell=1}^{L} \left\{\operatorname{diag}\left[S_\ell\right] W_{L-\ell}\right\}\biggr)\cdot h_1=\\
\biggl(3^{L} W_{L} 
\smm{i=1}{L}
\bigl(
\prod_{\ell=1}^{i-1} 
\left\{\operatorname{diag}\left[S_\ell\right] W_{L-\ell}\right\}
\bigr) d\left\{\operatorname{diag}\left[S_i\right] W_{L-i}\right\}
\cdot 
\bigl(
\prod_{\ell=i+1}^{L} \left\{\operatorname{diag}\left[S_\ell\right] W_{L-\ell}\right\}
\bigr)
\biggr)\cdot h_1\,.
\end{multline*}
Recall, that
$\operatorname{diag}(v) = \smm{i=1}{n}e_i^T ve_i e_i^T$, where $e_i$ - $i$-th standard basis vector.

\begin{align*}
d(\operatorname{diag}\left[S_{\ell}\right]) 
&= \smm{i=1}{p_{L-\ell+1}}e_i^T d((\AAA_{L-\ell,1}(x) + v_{L-\ell+1} \vee 0)^2)e_i e_i^T \\
&= \smm{i=1}{p_{L-\ell+1}}e_i^T 2\operatorname{diag}\left[\AAA_{L-\ell,1}(x) + v_{L-\ell+1} \vee 0\right]\cdot 
J(\AAA_{L-\ell,1}(x))\cdot 
dx \cdot e_i e_i^T\eqsp.
\end{align*}
Note that $J(\AAA_{L-\ell,1}(x)) = \nabla(\AAA_{L-\ell,1}(x))$ which we calculated above. By elementary algebra

\begin{align*}
&\langle d(\nabla f), h_1\rangle =
\biggl(3^{L} W_{L} 
\smm{i=1}{L}
\bigl(
\prod_{\ell=1}^{i-1} 
\left\{\operatorname{diag}\left[S_\ell\right] W_{L-\ell}\right\}
\bigr) 
d\left\{\operatorname{diag}\left[S_i\right] W_{L-i}\right\}
\bigl(
\prod_{\ell=i+1}^{L} \left\{\operatorname{diag}\left[S_\ell\right] W_{L-\ell}\right\}
\bigr)
\biggr)\cdot h_1=\\
&\biggl(3^{L} W_{L} 
\smm{i=1}{L}
\bigl(
\prod_{\ell=1}^{i-1} 
\left\{\operatorname{diag}\left[S_\ell\right] W_{L-\ell}\right\}
\bigr) \bigl(
\smm{j=1}{p_{L-i+1}}e_j^T 2\operatorname{diag}\left[\AAA_{L-i,1}(x) + v_{L-i+1} \vee 0\right]\cdot 
J(\AAA_{L-i,1}(x))\cdot 
dx \cdot e_j e_j^T
\bigr)
\cdot\\
&W_{L-i}
\bigl(
\prod_{\ell=i+1}^{L} \left\{\operatorname{diag}\left[S_\ell\right] W_{L-\ell}\right\}
\bigr)
\biggr)\cdot h_1\\
&=
(dx)^T
\biggl(
\smm{i=1}{L}2\smm{j=1}{p_{L-i+1}}
3^{L} \biggl[W_{L} \bigl(
\prod_{\ell=1}^{i-1} 
\left\{\operatorname{diag}\left[S_\ell\right] W_{L-\ell}\right\}
\bigr)
\cdot 
e_j^T \operatorname{diag}\left[\AAA_{L-i,1}(x) + v_{L-i+1} \vee 0\right]\cdot 
J(\AAA_{L-i,1}(x))\biggr]^T \cdot\\ 
&e_j e_j^T
\cdot
W_{L-i}
\bigl(
\prod_{\ell=i+1}^{L} \left\{\operatorname{diag}\left[S_\ell\right] W_{L-\ell}\right\}
\bigr)
\biggr)\cdot h_1=\langle \nabla^2 f \cdot dx, h_1\rangle\eqsp.
\end{align*}
From this expression Hessian is extracted.
\end{proof}

Now we need the counterpart of \Cref{Lem:inf_norm_bound} for the Jacobians $\nabla \fone(x)$, $\nabla \ftwo(x)$ and Hessians $\nabla^2 \fone(x)$, $\nabla^2 \ftwo(x)$. We prove the following lemma.
\begin{Lem}
\label{Lem:inf_norm_bound_derivative}
For the neural networks $\fone(x)$ and $\ftwo(x)$ defined in \eqref{eq:f_one_f_two_def}, $\|x\|_\infty \leq 1$ and output dimension 1, it holds
\begin{equation}
\label{eq:inf_norm_bound_2}
\norm{\nabla \fone(x) - \nabla \ftwo(x)}_{\infty} \leq \eps L (L+1) 3^{2L+2}
\bigl(
\prod_{\ell=0}^{L}(p_{\ell}+1)^{L3^{L+4}}
\bigr)\eqsp,
\end{equation}

\begin{equation}
\label{eq:inf_norm_bound_hes_2}
\norm{\nabla^2 \fone(x) - \nabla^2 \ftwo(x)}_{\infty} \leq \eps L^2(L+1)(L+2)3^{3L+4}
\bigl(
\prod_{\ell=0}^L
(p_\ell+1)^{L3^{L+5}}
\bigr)\,.
\end{equation}
\end{Lem}
\begin{proof}
Recall the expression of derivative \eqref{eq:first_grad}.
Let us define for $j \in \{0,\dots,L-1\}$ and $i \in \{1,2\}$ the quantities 
\begin{align*}
\Delta^{(i)}_{j} &= (\AAA_{j,1}^{(i)}(x) + v^{(i)}_{j+1} \vee 0)^2\,.
\end{align*}
Note that $\diagone_{j}, \diagtwo_{j} \in \R^{p_{j+1}}$. With the triangular inequality,
\begin{align*}
&\norm{\nabla \fone(x) - \nabla \ftwo(x)}_{\infty} 
\leq 3^{L} \norm{\bigl(\Wone_{L} - \Wtwo_{L}\bigr)\prod_{\ell=1}^{L}\left\{\operatorname{diag}\left[\diagone_{L-\ell}\right] \Wone_{L-\ell}\right\}}_{\infty} 
\\&+ 3^{L}(p_{L}+1) \times 
\sum_{k=1}^{L} \stnorm{\prod_{\ell=1}^{k-1}\left\{\operatorname{diag}\left[\diagone_{L-\ell}\right] \Wone_{L-\ell}\right\}\\
&\times \left(\operatorname{diag}\left[\diagone_{L-k}\right] \Wone_{L-k} - \operatorname{diag}\left[\diagtwo_{L-k}\right] \Wtwo_{L-k}\right)
\prod_{\ell=k+1}^{L}\left\{\operatorname{diag}\left[\diagtwo_{L-\ell}\right] \Wtwo_{L-\ell}\right\}}_{\infty}.
\end{align*}
Proceeding as in \Cref{Lem:inf_norm_bound}, we obtain 
\begin{align*}
&\norm{\operatorname{diag}[\diagone_{k}] \Wone_{k} - \operatorname{diag}[\diagtwo_{k}] \Wtwo_{k}}_{\infty} \leq 
\norm{\operatorname{diag}[\diagone_{k} - \diagtwo_{k}]\Wtwo_{k}}_{\infty} 
+
\norm{\operatorname{diag}[\diagone_{k}] (\Wone_{k} -\Wtwo_{k})}_{\infty} \\
&\qquad\leq 
\bigl(\eps(k+1)3^{k}\prod_{\ell=0}^{k}(p_{\ell}+1)^{3^{k+2}}\bigr)
(p_k+1)
\cdot 2(p_k+1)
\bigl(
\prod_{\ell=1}^{k}(p_{\ell}+1)^{3^{k}}
\bigr)\\
&\qquad + \eps(p_k+1)(p_k+1)^2
\bigl(
\prod_{\ell=1}^{k}(p_{\ell}+1)^{3^{k}}
\bigr)^2
\leq 2\eps (k+1)3^{k+1} \prod_{\ell=0}^{k}(p_{\ell}+1)^{3^{k+3}}\,,
\end{align*}
and, similarly,
\begin{equation*}
\norm{\operatorname{diag}\left[\diagone_{L-\ell}\right] \Wone_{L-\ell}}_{\infty} \\
\leq (p_k+1)\cdot \bigl((p_k+1)\prod_{j=1}^{L-\ell}(p_{L-j}+1)^{3^{j}}\bigr)^2
\leq \prod_{j=1}^{L-\ell}(p_{L-j}+1)^{3^{j+1}}\,.
\end{equation*}
Combining the previous bounds, we obtain
\begin{align*}
\norm{\nabla \fone(x) - \nabla \ftwo(x)}_{\infty} 
&\leq 3^L (p_L+1) L\cdot 2\cdot \eps (L+1)3^{L+1}
\bigl(
\prod_{\ell=1}^{L}(p_{\ell}+1)^{3^{L+3}}
\bigr)
\bigl(
\prod_{\ell=0}^{L}(p_{\ell}+1)^{L3^{L+1}}
\bigr) \\
&+\eps 3^{L}\prod_{\ell=0}^{L}(p_{\ell}+1)^{L3^{L+1}} \leq 
\eps L (L+1) 3^{2L+2}
\bigl(
\prod_{\ell=0}^{L}(p_{\ell}+1)^{L3^{L+4}}
\bigr)
\,,  
\end{align*}
and the statement follows. The same is done for hessian matrices, thus some steps are omitted.
\begin{align*}
\stnorm{\nabla^2 \fone(x) - \nabla^2 & \ftwo(x)}_{\infty} 
\leq 2\cdot 3^L L\cdot (\max_{i\leq L} p_i)\cdot (2L+4) \bigl(
\prod_{\ell=0}^L
\bigl(
(p_\ell+1)
\prod_{s=1}^{L-\ell}
(p_{L-\ell-s}+1)^{3^s}
\bigr)^5
\bigr)
\biggr) \\
&\times
\max
\biggl(
\eps (\max_{i\leq L} p_i),
\eps (L+1) 3^{L+1} \prod_{\ell=0}^{L}
(p_{\ell}+1)^{3^{L+3}}, \eps L(L+1) 3^{2L+2} \prod_{\ell=0}^{L}
(p_{\ell}+1)^{L3^{L+4}} \biggr)\\
&\leq 
\eps L^2(L+1)(L+2)3^{3L+4}
\bigl(
\prod_{\ell=0}^L
(p_\ell+1)^{L3^{L+5}}
\bigr).
\end{align*}
\end{proof}

\section{Details of numerical experiments}
\label{sec:numerical_exp_details}	
	\begin{table}[H]
        \centering
        \begin{tabular}{ |c|c|c|c|c|c|c| } 
        \hline
        CGM & 
        $\gamma$&
        $n_{\text{burn}}$& 
        $n_{\text{train}}$&
        $n_{\text{test}}$&
        $T$&
        $b_n$\\
        \hline
        NUTS & 
        0.1&
        $10^4$& 
        $3\cdot 10^4$&
        $3\cdot 10^4$&
        30&
        30\\
        \hline
        \end{tabular}
        \caption{Funnel, algorithm and chain generation hyperparameters}
    \end{table}
    \begin{table}[H]
    \centering
        \begin{tabular}{ |c|c|c| } 
        \hline
        $d$ & $a$ & 
        $b$\\
        \hline
        2& 1&
        0.5\\
        \hline
        \end{tabular}
        \caption{Funnel, distribution hyperparameters}
    \end{table}
    \begin{table}[H]
    \centering
        \begin{tabular}{ |c|c|c|c|c|c|c| } 
        \hline
        CGM & 
        $\gamma$&
        $n_{\text{burn}}$& 
        $n_{\text{train}}$&
        $n_{\text{test}}$&
        $T$&
        $b_n$\\
        \hline
        ULA & 
        0.01&
        $10^5$& 
        $2\cdot 10^4$&
        $10^4$&
        30&
        30\\
        \hline
        \end{tabular}
        \caption{Banana-shaped, algorithm and chain generation hyperparameters}
    \end{table}
    \begin{table}[H]
    \centering
        \begin{tabular}{ |c|c|c|c|c| } 
        \hline
        $d$ & 
        $p$&
        $b$\\
        \hline
        6 & 
        20 & 
        $0.05$\\
        \hline
        \end{tabular}
        \caption{Banana-shaped, distribution hyperparameters}
    \end{table}
    \begin{table}[H]
    \centering
        \begin{tabular}{ |c|c|c|c|c|c|c| } 
        \hline
        CGM & 
        $\gamma$&
        $n_{\text{burn}}$& 
        $n_{\text{train}}$&
        $n_{\text{test}}$&
        $T$&
        $b_n$\\
        \hline
        ULA & 
        0.1&
        $10^4$& 
        $3\cdot10^4$&
        $10^4$&
        30&
        15\\
        \hline
        \end{tabular}
        \caption{Logistic Regression, algorithm and chain generation hyperparameters}
    \end{table}
    \begin{table}[H]
    \centering
        \begin{tabular}{ |c|c|c| } 
        \hline
        $d$ & 
        $g$&
        $K$\\
        \hline
        9 & 
        100 & 
        154\\
        \hline
        \end{tabular}
        \caption{Logistic Regression, distribution hyperparameters}
    \end{table}

\newpage
\bibliographystyle{elsarticle-harv}
\bibliography{bibfile-mr}       

\end{document}